\documentclass[onecolumn,rmp,11pt]{revtex4}

\usepackage{graphicx, amsmath, amssymb, bm, enumerate, mathtools}

\usepackage[hyperfootnotes=false]{hyperref}
\usepackage{xcolor}
\hypersetup{
    colorlinks,
    linkcolor={green!60!black},
    citecolor={blue!90!black},
}
\begin{document}

\title{Dimension reduction for stochastic dynamical systems forced onto a manifold by large drift: a constructive approach with examples from theoretical biology}

\author{Todd L. Parsons}
\affiliation{Laboratoire de Probabilit\'{e}s et Mod\`{e}les Al\'{e}atoires, CNRS UMR 7599, Universit\'{e} Pierre et Marie Curie, Paris, 75005, France.}
\author{Tim Rogers}
\affiliation{Centre for Networks and Collective Behaviour, Department of Mathematical Sciences, University of Bath, Bath, BA2 7AY, UK.}
\pacs{}
\begin{abstract}
Systems composed of large numbers of interacting agents often admit an effective coarse-grained description in terms of a multidimensional stochastic dynamical system, driven by small-amplitude intrinsic noise. In applications to biological, ecological, chemical and social dynamics it is common for these models to posses quantities that are approximately conserved on short timescales, in which case system trajectories are observed to remain close to some lower-dimensional subspace. Here, we derive explicit and general formulae for a reduced-dimension description of such processes that is exact in the limit of small noise and well-separated slow and fast dynamics. The Michaelis-Menten law of enzyme-catalyzed reactions, and the link between the Lotka-Voltera and Wright-Fisher processes are explored as a simple worked examples. Extensions of the method are presented for infinite dimensional systems and processes coupled to non-Gaussian noise sources.
\end{abstract}
\maketitle

\section{Introduction}

One important reason for the observed ``unreasonable effectiveness'' of mathematical modelling in describing natural phenomena \cite{Wigner1960} is the gigantic separation of scales apparent in physical systems. This convenient property allows a process of interest to be treated separately from what is happening on much slower or faster time scales, or on much larger or smaller spatial scales. The result is that mathematical models can be simple to state and, having few degrees of freedom, are more likely to be solvable. Sadly, other fields to which mathematics is applied --- biology, ecology, finance, social science, to name a few --- do not enjoy this separation, and consequently modellers in these areas have not yet come close to the same levels of predictive success achieved in physics. In writing a successful mathematical model, one faces two intertwined problems: the ability of a model to capture the essential behaviours of the system of interest, and our ability to `solve' the model to extract quantitative predictions and qualitative understanding. Biological systems, for example, often comprise feedbacks over several scales, meaning that more degrees of freedom must be included in a model to achieve an acceptable representation of reality. This has negative consequences for the analysis, however, as higher-dimensional systems tend to be harder to work with, both analytically and in simulations.

These concerns are particularly relevant for the study of emergent phenomena arising from the interactions of large numbers of individual elements (typically particles, agents, or organisms). Important examples include the modelling of population dynamics, evolution, and epidemic spread. Significant progress has been made in this area through the application of so-called system size expansion techniques (e.g. \cite{vanKampen92,Kurtz78}). These methods provide a coarse-grained description of an interacting system, typically in terms of a stochastic process driven by intrinsic noise (i.e. with amplitude depending on the state of the system, and scaling inversely with system size), which is precise in the limit of large system size. Although very much less complicated than the original system, these coarse-grained models are still likely to resist exact solution, being non-linear multi-dimensional stochastic dynamical systems. Previously, important results have been obtained though analysis of linear perturbations around fixed points, however, to more fully characterise the phenomena possible in these models requires a fully non-linear treatment.

\begin{figure}
\includegraphics[width=0.8\textwidth, trim=50 280 50 280]{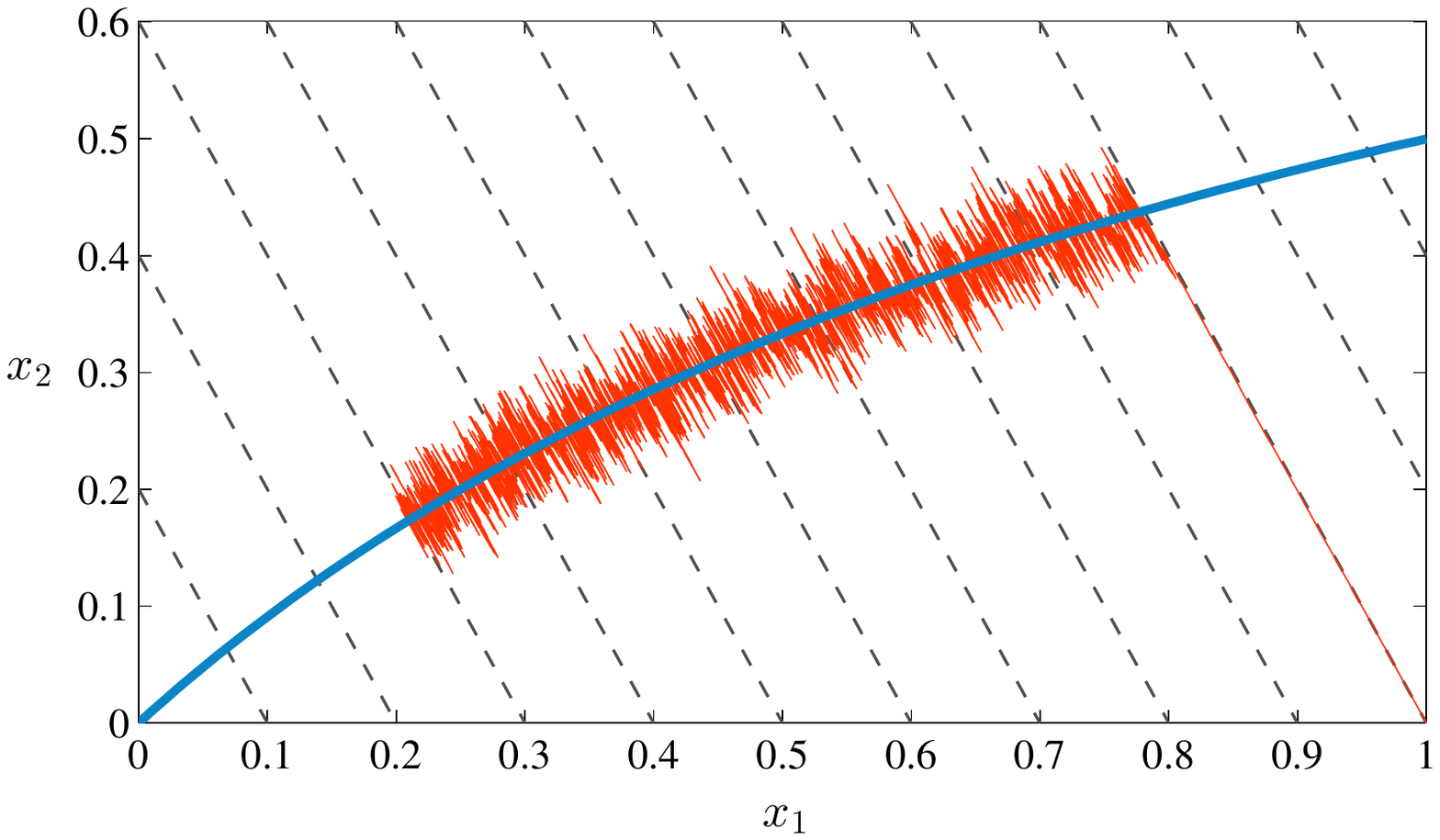}
\caption{\textbf{Thin Red:} Simulation of a single stochastic trajectory of an SDE of the type (\ref{sdex}), with $\bm{f}$, $\bm{h}$ and $\bm{G}$ corresponding to the Michaelis-Menten model (\ref{MMfh}, \ref{MMG}). \textbf{Thick Blue:} The slow manifold for this system, which the stochastic trajectory stays close to after the fast initial transient carrying it away from the initial condition $(1,0)$. \textbf{Dashed Black:} The flow field of the outer drift term $\bm{f}$, to which the fast motion is approximately parallel.}
\label{MMfig}
\end{figure}

Working in this direction, several groups have independently found and exploited a natural separation of scales emerging in certain models of interacting populations \cite{Parsons2007b,Parsons2008,Durrett2009,Parsons2010,Parsons2012,Rogers2012a,Rogers2012b,Doering2012,Constable2013,Constable2014,Pigolotti2014,Constable2015,Chotibut2015}. Loosely speaking, it is often the case that the total size of a population varies much more rapidly than its composition, as is evident in the disparate timescales of ecology and evolution. In dynamical systems terminology, it is observed that trajectories remain in the neighbourhood of a lower-dimensonal manifold; a subspace of the system state space in which the total size of the population is a function of its composition. Intrinsic noise drives small perturbations from this manifold, which are quickly suppressed by a large deterministic drift back (see, for example, the trajectories of Michaelis-Menten dynamics in Fig.~\ref{MMfig}). The works cited above pursue various related approximation strategies, allowing for a simplified, often solvable, effective model to be derived describing motion along the lower-dimensional manifold.

This is a kind of timescale separation that cannot be put by hand into a model as there is a complex feedback between the fast and slow degrees of freedom, which must be carefully computed. Mathematically, this is not straightforward, and a host of different strategies have previously been employed to approximate the effects of this feedback between scales. A rigorous treatment of this situation can in fact be found in the theorems proved in \cite{Katzenberger91} and \cite{Funaki1993}, however, the results of these works are difficult to apply in practice as they lack explicit formulations of certain key quantities. Our purpose here is to synthesise the various methods of the above mentioned authors with the rigorous theory of Katzenberger, and distil this into a single robust, systematic and provably correct procedure for timescale separation in stochastic dynamical systems with intrinsic noise, which we believe will be of considerable general use. The main results are contained in the short Section \ref{S2}, where we describe a map from a high-dimensional system of equations (\ref{sdex}) to a lower-dimensional one (\ref{sdez}), via explicit formulae that are summarised in Table~\ref{tab1}.
\begin{table}
\begin{center}
\begin{tabular}{l|l}
\textbf{Case} & \textbf{Procedure} \\\hline\hline
Outer system ($\varepsilon=\mu=0$) is solvable $\quad$& Use equations (\ref{ivp}, \ref{defpi}) and (\ref{defPQ}, \ref{drift})\\
Manifold is one-dimensional & Use equations (\ref{drift}) and (\ref{P1D}, \ref{Q1D})\\
Manifold has co-dimensional one & Use equations (\ref{Pcd1}) and (\ref{Qcd1})\\
Manifold is $m$-dimensional & Use equations (\ref{drift}) and (\ref{diagJ} -- \ref{QmD})
\end{tabular}
\end{center}
\caption{Quick reference table of equations applying to different cases of slow-manifold reduction.}
\label{tab1}
\end{table}

The general setting for our calculations will be models expressed as coupled stochastic differential equations, with some small parameter $\varepsilon$ that controls the separation of timescales. As stated above, we are particularly interested in systems that are derived from a complex underlying interacting process, so that the noise terms are intrinsic, representing the cumulative random effects of many interactions. Importantly in this case the noise, although non-negligible in its effect, is typically small in amplitude\footnote{A complementary branch of theory exists dealing with the relaxation of this assumption, see \cite{Arnold1998,Roberts2008} for starting points in the literature.}. Three example application areas are: (i) chemical and biological reaction networks, \textit{e.g.}  gene regulation in a cell t, where small copy numbers imply noise but homoeostasis suggests timescale separation \cite{Ball2006}, (ii) evolutionary models where noise arising from demographic fluctuations can alter the course of selection (references above), (iii) dynamical networks which are naturally extremely high-dimensional systems in need of low-dimensional proxies \cite{Rogers2013}. It is worth pointing out that a plethora of different timescale separation techniques exist in the literature, and the most useful choice of method depends on the system in question. In Appendix \ref{lit} we give a very brief history of notable developments in stochastic timescale separation.

In addition to the general theory outlined in Section~\ref{S2}, we present in Section~\ref{WE} some worked examples and variations to the method. As a prototypical example of the standard method, we work through the derivation of a stochastic form of the Michaelis-Menten law for enzyme-catalysed reactions. The second example demonstrates the extension of the methods to an infinite dimensional setting. Katzenburger's theorem is discussed in detail in Appendix \ref{Katz}, as well as the generalisation of the noise sources (for example to jump processes) and alternative characterisations of the diffusion process on the manifold.

\section{Reduced model description}\label{S2}

We consider Langevin stochastic differential equations of the general type
\begin{equation}
\frac{d\bm{x}}{dt}=\bm{f}(\bm{x})+\varepsilon\bm{h}(\bm{x})+\sqrt{\mu}\,\bm{G}(\bm{x})\,\bm{\eta}(t)\,,
\label{sdex}
\end{equation}
where the state variable $\bm{x}=(x_1,\ldots,x_d)^T$ is an $d$-dimensional vector, and there are $s$ independent It\^{o} white noise sources $\bm{\eta}(t)=(\eta_1(t)\,,\ldots,\eta_s(t))^T$. The vector-valued functions $\bm{f}:\mathbb{R}^d\to\mathbb{R}^d$ and $\bm{h}:\mathbb{R}^d\to\mathbb{R}^d$ are the `outer' and `inner' parts of the drift respectively, and the matrix valued function $\bm{G}:\mathbb{R}^d\to\mathbb{R}^{d,s}$ specifies the coupling of state variables to noise sources. We assume throughout that $\bm{f}$ is twice differentiable, but place no constraints on the other functions. The parameters $\varepsilon$ and $\mu$ determine the separation of timescales and the strength of the noise, respectively.

We do not assume an \textit{a priori} separation into slow and fast variables, as is common in the literature, as in the applications that motivate us, an appropriate change of variables is frequently neither evident nor analytically tractable (although see \cite{Parsons2007b}), and our method does not require that they be known.

We are interested in the case when $\varepsilon$ and $\mu$ are small and $\bm{f}$ possesses a $m$-dimensional manifold of equilibria $\Gamma\subset\mathbb{R}^d$ such that $\bm{f}(\bm{\tilde{x}})=\bm{0}$ for all $\bm{\tilde{x}}\in\Gamma$.  We assume that this manifold is unique, connected, and globally attractive (\textit{i.e.} it is a globally unique, normally hyperbolic slow manifold see \textit{e.g} \cite{Berglund2006}); then we expect solutions of (\ref{sdex}) to rapidly approach and remain very close to $\Gamma$. In fact, it has been rigourously proved by Katzenberger \cite{Katzenberger91} that the trajectories of $\bm{x}\in\mathbb{R}^d$ converge those of a stochastic variable $\bm{\tilde{x}}\in\Gamma$ with dynamics
\begin{equation}
\frac{d\bm{\tilde{x}}}{dt}=\varepsilon\bm{P}(\bm{\tilde{x}})\bm{h}(\bm{\tilde{x}})+\mu\bm{g}(\bm{\tilde{x}})+ \sqrt{\mu}\bm{P}(\bm{\tilde{x}})\bm{G}(\bm{\tilde{x}})\bm{\eta}(t)\,,
\label{sdez}
\end{equation}
where $\bm{P}$ is a certain projection matrix derived from $\bm{f}$, and $\bm{g}$ is a new contribution to the drift arising from the way in which fluctuations away from the manifold are suppressed; as our examples illustrate, unlike the deterministic situation, it is not sufficient to simply restrict \eqref{sdex} to $\Gamma$ to obtain the slow dynamics. Our purpose here is to derive explicit expressions for $\bm{P}$ and $\bm{g}$. Readers with a specific problem in mind may wish to jump straight to the appropriate result, which can be found by referring to Table~\ref{tab1}.

Before we proceed with our main task, we give a brief sketch of the derivation of (\ref{sdez}). Examining (\ref{sdex}) when $\varepsilon$ and $\mu$ are small, one might imagine a picture in which the state of the system is quickly carried onto the manifold by the fast outer drift term $\bm{f}$. Following this fast initial transient, it may then receive multiple stochastic `kicks' carrying it away from the manifold, each time only to return again via the paths described by $\bm{f}$. See Figure \ref{MMfig} for an illustrative example. This intuition can be made concrete by considering the flow map of the outer system. Let $\bm{x}$ be a point in the state space and consider the deterministic initial value problem
\begin{equation}
\begin{cases}\displaystyle\frac{d\bm{\xi}_{\bm{x}}}{dt}=\bm{f}(\bm{\xi}_{\bm{x}})\\ \bm{\xi}_{\bm{x}}(0)=\bm{x}\,.\end{cases}
\label{ivp}
\end{equation}
Since the center manifold is globally attractive, all trajectories lead eventually to $\Gamma$ and we may thus define a function $\bm{\pi}:\mathbb{R}^d\to\Gamma$ giving the endpoint of the deterministic trajectories
\begin{equation}
\bm{\pi}(\bm{x})=\lim_{t\to\infty}\bm{\xi}_{\bm{x}}(t)\,,
\label{defpi}
\end{equation}
where $\bm{\xi}_{\bm{x}}$ is the solution of (\ref{ivp}).

If we take the point $\bm{x}$ to be the current location of the random variable governed by equation (\ref{sdex}), then $\bm{\pi}(\bm{x})$ defines another random variable that tracks the motion of $\bm{x}$ but is constrained to the manifold. Application of It\^{o}'s formula \cite{Ito1974} gives the Langevin equations for each spatial coordinate:
\begin{equation}
\begin{split}
\frac{d}{dt}\pi_i(\bm{x})&=\sum_j \frac{\partial\pi_i}{\partial x_j} \frac{dx_j}{dt} + \frac{\mu}{2}\sum_{s,j,k} G_{js}(\bm{x})G_{ks}(\bm{x})\frac{\partial^2\pi_i}{\partial x_j\partial x_k}\\
&=\sum_j \frac{\partial\pi_i}{\partial x_j}\Big(f_j(\bm{x})+\varepsilon h_j(\bm{x})\Big) +\frac{\mu}{2}\sum_{s,j,k} G_{js}(\bm{x})G_{ks}(\bm{x})\frac{\partial^2\pi_i}{\partial x_j\partial x_k} +\sqrt{\mu}   \sum_{s,j}G_{js}(\bm{x})\frac{\partial\pi_i}{\partial x_j}\eta_s(t)\\&
=\varepsilon\sum_j \frac{\partial\pi_i}{\partial x_j} h_j(\bm{x}) +\frac{\mu}{2}\sum_{s,j,k} G_{js}(\bm{x})G_{ks}(\bm{x})\frac{\partial^2\pi_i}{\partial x_j\partial x_k} +\sqrt{\mu}   \sum_{s,j}G_{js}(\bm{x})\frac{\partial\pi_i}{\partial x_j}\eta_s(t)
\,.
\label{sdezx}
\end{split}
\end{equation}
where the last equality comes from the observation that $\bm{\pi}(\bm{\xi}_{\bm{x}}(t)) = \bm{\pi}(\bm{x})$ for all $t$, and thus
\begin{equation}
	0 = \frac{d}{dt}\Big\vert_{t = 0} \pi_{i}(\bm{\xi}_{\bm{x}}(t))
		= \sum_j \frac{\partial\pi_i}{\partial x_j} \frac{d\xi_{j}}{dt}\Big\vert_{t = 0}
		= \sum_j \frac{\partial\pi_i}{\partial x_j}  f_j(\bm{x})\,.
\end{equation}
Unfortunately equation \eqref{sdezx} is not closed since it relies on full knowledge of the random variable $\bm{x}$. However, if we believe that $\bm{x}$ remains very close to $\Gamma$ (as is the case when $\varepsilon$ and $\mu$ are small) then we might be motivated to consider a new random variable $\bm{\tilde{x}}\in\Gamma$ which we assume is a close approximation to both $\bm{x}$ and $\bm{\pi}(\bm{x})$. Substituting $\bm{\tilde{x}}$ for both these quantities in \eqref{sdezx}, we obtain the closed expression
\begin{equation}
\frac{dz_i}{dt}=\varepsilon \sum_j P_{ij}(\bm{\tilde{x}}) h_j(\bm{\tilde{x}}) +\frac{\mu}{2} \sum_{s,j,k}G_{js}(\bm{\tilde{x}})G_{ks}(\bm{\tilde{x}})Q_{ijk}(\bm{\tilde{x}}) +\sqrt{\mu} \sum_{s,j}P_{ij}(\bm{\tilde{x}})G_{js}(\bm{\tilde{x}})\eta_s(t)\,,
\label{sdez2}
\end{equation}
where $\bm{P}$ is a matrix and $\bm{Q}$ an array defined by
\begin{equation}
P_{ij}(\bm{\tilde{x}})=\left.\frac{\partial}{\partial x_j}\pi_i(\bm{x})\right|_{\bm{x}=\bm{\tilde{x}}}\,,\quad Q_{ijk}(\bm{\tilde{x}})=\left.\frac{\partial^2}{\partial x_j\partial x_k}\pi_i(\bm{x})\right|_{\bm{x}=\bm{\tilde{x}}}\,.
\label{defPQ}
\end{equation}
Equivalently we may rewrite (\ref{sdez2}) as equation (\ref{sdez}), where the additional drift term is
\begin{equation}\label{drift}
\bm{g}(\bm{\tilde{x}})=\frac{1}{2}\sum_{s,j,k} G_{js}(\bm{\tilde{x}})G_{ks}(\bm{\tilde{x}})Q_{ijk}(\bm{\tilde{x}}).
\end{equation}

\begin{figure}
\includegraphics[width=\textwidth, trim=0 0 0 0]{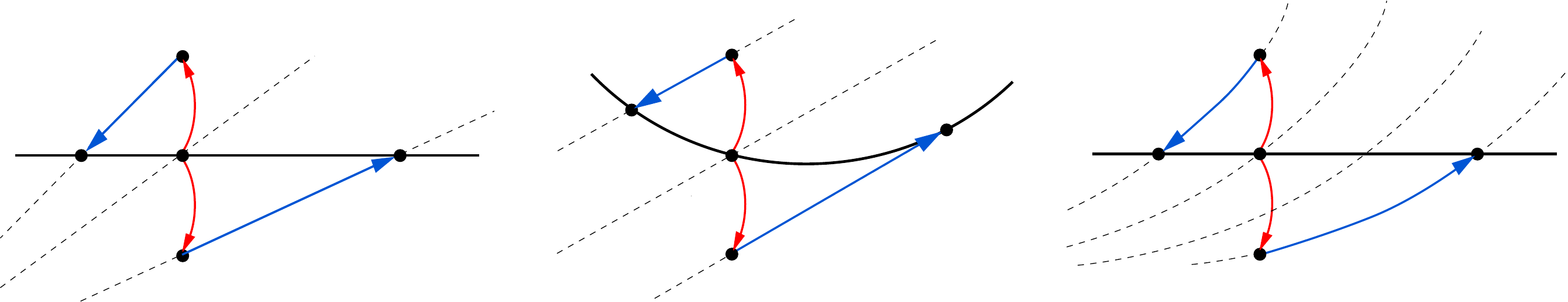}
\caption{\textbf{Left:} Here, the variation in the angle between the fast (dashed) and slow (solid) subspaces creates a bias in the location of the return to the manifold of a perturbation away from it; an upward perturbation returns quite close on the left of the origin, but an equally likely downward perturbation is carried far to the right. \textbf{Centre:} The same effect can occur as a result of curvature of the manifold. In this figure the flow fields are parallel, but the manifold curves, resulting in the same rightwards bias in the projected system. \textbf{Right:} Curvature of the flow field may also induce bias, even when the angle of intersection is constant.}
\label{ffd_fig}
\end{figure}
The projection matrix $\bm{P}(\bm{\tilde{x}})$ is entirely determined by the first order terms of this expansion, and typically it can be straightforwardly reconstructed from knowledge of the eigenvectors of the Jacobian matrix of $\bm{f}$. The calculation of the noise-induced drift term is more complicated, having contributions from three possible sources: variation of the alignment of the flow field, curvature of the manifold, and curvature of the flow field. Each of these mechanisms can induce a bias in the direction of flow of the reduced dimension system, as illustrated in Figure \ref{ffd_fig}.

In the following subsections we will present explicit procedures for computing $\bm{P}$ and $\bm{Q}$.

\subsection{One-dimensional manifolds}

The simplest case to treat is that of a one-dimensional manifold, as the second-order perturbation expansion is explicitly solvable. Suppose that the slow manifold $\Gamma$ is a curve parameterised by the first spatial co-ordinate of the system\footnote{Note that we have chosen this case for simplicity of presentation, and not all 1D manifolds can in fact be tackled in this way (e.g. a circular manifold would fail here). The more general case of a manifold described by an arbitrary parameterised curve is not substantially different, however, as we only ever require the local properties of the projection $\pi$, and for smooth manifolds there is always a local coordinate system in which the problem can be set up in the required format.}.  That is, there exists function $\bm{\gamma}$ such that
\begin{equation}
\bm{x}\in\Gamma\quad\Leftrightarrow\quad \bm{x}=\bm{\gamma}(x_1)\,.
\end{equation}
In this case the dynamics of the reduced system $\bm{\tilde{x}}$ defined in (\ref{sdez}) are determined entirely by the first component, so we need only to compute the partial derivatives of $\pi_1$. For ease of notation we will drop the subscript $1$ from now on, writing $\tilde x:=\tilde x_1$ as well as $P_{j}:=P_{1j}(\bm{\gamma}(\tilde x_1))$ and $Q_{jk}:=Q_{1jk}(\bm{\gamma}(\tilde x_1))$.

\begin{figure}
\includegraphics[width=0.8\textwidth, trim=0 0 0 0]{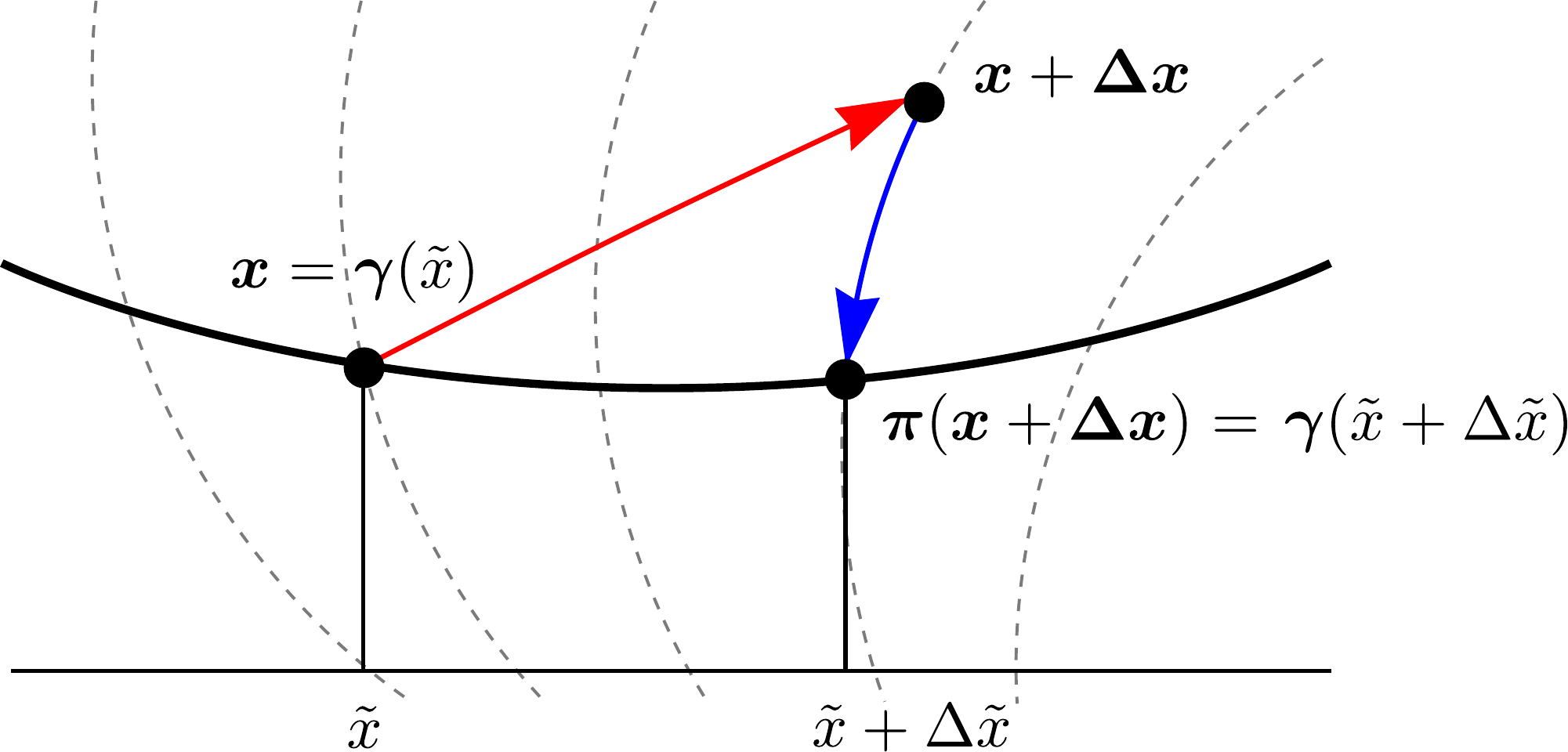}
\caption{Illustration of the perturbation calculation for a 1D manifold. }
\label{1Dfig}
\end{figure}

Consider a point $\bm{x}=\bm{\gamma}(\tilde x)$ on the manifold. To obtain expressions for $P_{j}$ and $Q_{jk}$ we undertake a second-order perturbation theory. Since $\bm{x}\in\Gamma$ is a point on the manifold we have by definition that it is unchanged by the action of the outer flow field, so $\bm{\pi}(\bm{x})=\bm{\gamma}(\tilde x)$. We make a small perturbation $\bm{x}\mapsto \bm{x}+\bm{\Delta x}$, and ask what perturbation $\tilde x\mapsto \tilde x+\Delta \tilde x$ is required so that $\bm{\pi}(\bm{x}+\bm{\Delta x})=\bm{\gamma}(\tilde x+\Delta \tilde x)$. See Figure~\ref{1Dfig} for an illustration.

Becuase we set the problem up so that $\pi_1(\bm{x})=\tilde x$, Taylor expanson gives
\begin{equation}
\Delta \tilde x=\sum_j P_{j} \Delta x_j+\frac{1}{2}\sum_{j,k}Q_{jk}\Delta x_j\Delta x_k+\ldots
\end{equation}
Near the point $\bm{x}\in \Gamma$ we can approximate the action of $\bm{\pi}$ by constructing the quadratic expansion of the preimage. Specifically, it can be shown that in the neighbourhood of $\bm{x}$, the collection of nearby points that would be mapped to $\bm{x}$ by $\bm{\pi}$ (i.e. the set $\{\bm{y}\,:\,\bm{\pi}(\bm{y})=\bm{x}\}$) is approximated to second order by the set of points $\bm{y}$ such that
\begin{equation}
\bm{v}(\tilde x)^T(\bm{y}-\bm{\gamma}(\tilde x))+(\bm{y}-\bm{\gamma}(\tilde x))^T\Theta(\tilde x)(\bm{y}-\bm{\gamma}(\tilde x))=0\,,
\label{quadratic}
\end{equation}
where $\bm{v}(\tilde{x})$ is a perpendicular vector to the flow field near $\bm{x}=\bm{\gamma}(\tilde x)$ and $\Theta(\tilde{x})$ a matrix describing the curvature of the flow field near the same point. In Appendix \ref{Theta} we give an explicit derivation of these quantities from $\bm{f}$; for now we assume they are known. Recall that we are seeking the perturbation $\Delta\tilde{x}$ such that $\bm{\pi}(\bm{x}+\bm{\Delta x})=\bm{\gamma}(\tilde x+\Delta \tilde x)$, to second order. We make the following Taylor expansions of various orders:
\begin{equation}
\begin{split}
\Big[\bm{x}+\Delta\bm{x}-\bm{\gamma}(\tilde x+\Delta \tilde x)\Big]_\ell &= \Delta x_\ell- \gamma_\ell'\Delta\tilde x-\frac{1}{2}\gamma_\ell''(\Delta \tilde x)^2 +\ldots\\
&= \sum_k (\delta_{k,\ell}- \gamma_\ell'P_{k}) \Delta x_k -\frac{1}{2}\sum_{j,k} \left(\gamma_\ell' Q_{jk} + \gamma_\ell'' P_{j}P_k\right)\Delta x_j\Delta x_k+\mathcal{O}(\bm{\Delta x}^3)\\
\Big[\bm{v}(\tilde x+\Delta \tilde x)\Big]_\ell &=v_\ell+v'_\ell\Delta\tilde x+\ldots =v_\ell+v'_\ell\sum_j P_{j} \Delta x_j+\mathcal{O}(\bm{\Delta x}^2)\\
\Big[\Theta(\tilde x+\Delta \tilde x)\Big]_{jk} &= \Theta_{jk}+\mathcal{O}(\bm{\Delta x})\,.
\end{split}
\end{equation}
Here we drop the argument $\tilde x$ from $\gamma_\ell'$, $\gamma_\ell''$, $v_\ell$, $v_\ell'$ and $\Theta_{jk}$ to avoid clutter. Following (\ref{quadratic}), the requirement that $\bm{\pi}(\bm{x}+\bm{\Delta x})=\bm{\gamma}(\tilde x+\Delta \tilde x)$ to second order becomes
\begin{equation}
\begin{split}
0&=\bm{v}(\tilde x+\Delta \tilde x)^T\Big(\bm{x}+\bm{\Delta x}-\bm{\gamma}(\tilde x+\Delta \tilde x)\Big)\\
&\quad+\Big(\bm{x}+\bm{\Delta x}-\bm{\gamma}(\tilde x+\Delta \tilde x)\Big)^T\Theta(\tilde x+\Delta \tilde x)\Big(\bm{x}+\bm{\Delta x}-\bm{\gamma}(\tilde x+\Delta \tilde x)\Big)+\mathcal{O}(\bm{\Delta x}^3)\\
&=\sum_\ell \Big[v_\ell+v'_\ell\sum_j P_{j} \Delta x_j\Big]\Big[\sum_k (\delta_{k,\ell}- \gamma_\ell'P_{k}) \Delta x_k -\frac{1}{2}\sum_{j,k} \left(\gamma_\ell' Q_{jk} + \gamma_\ell'' P_{j}P_k\right)\Delta x_j\Delta x_k\Big]\\
&\quad+\sum_{j,k}\Theta_{jk}\Delta x_j\Delta x_k+\mathcal{O}(\bm{\Delta x}^3)\\
&=\sum_k\Big\{\sum_\ell   v_\ell(\delta_{k,\ell}- \gamma_\ell'P_{k}) \Big\}\Delta x_k\\
&\quad +\frac{1}{2}\sum_{j,k} \Big\{ \sum_\ell v_\ell'   (\delta_{k,\ell}+\delta_{j,\ell}- 2\gamma_\ell'P_{k})P_{j}-\sum_\ell v_\ell\left(\gamma_\ell' Q_{jk} + \gamma_\ell'' P_{j}P_k\right)+2\Theta_{jk}\Big\}\Delta x_j\Delta x_k \\&\qquad+\mathcal{O}(\bm{\Delta x}^3)\,.
\end{split}
\end{equation}
Since the perturbation $\bm{\Delta x}$ was abritrary, we require each term in curly brackets above to be equal to zero. From the first order terms we conclude that
\begin{equation}\label{P1D}
P_k=\frac{v_k}{\sum_\ell v_\ell\gamma_\ell'}\,,
\end{equation}
and from the second order that
\begin{equation}\label{Q1D}
Q_{jk}=\frac{1}{\sum_\ell v_\ell\gamma_\ell'}\left( v_k'P_j+v_j'P_k+2\Theta_{jk}  - \sum_\ell (2v_\ell' \gamma_\ell'+v_\ell\gamma_\ell'' ) P_{j}P_k\right).
\end{equation}
Written this way, the separate contributions from variation of the flow field (terms involving $\bm{v'}$), curvature of flow field ($\Theta$), and curvature of the manifold (the $\bm{\gamma''}$ term) are clearly visible.

In higher dimensions, the above perturbation expansion is less useful, as it produces a larger system of equations which lacks an explicit solution. A different line of attack is necessary.

\subsection{General case}
If the linearisation of the flow field $\bm{\phi}_t$ is known in the neighbourhood of the manifold then $\bm{P}$ can be reconstructed easily. Specifically, around a point $\bm{z}\in\Gamma$ the state space $\mathbb{R}^d$ can be decomposed into a product of `slow' and `fast' subspaces of dimension $m$ and $d-m$, respectively. The slow subspace is the tangent plane to the manifold at the given point; a perturbation in one of these directions is unaffected by the action of $\bm{f}$. Conversely, the fast subspace comprises perturbation directions that collapse quickly back to the manifold. The projection matrix $\bm{P}(\bm{z})$ acts as the identity on the slow subspace and as zero on the fast subspace.

Unfortunately, no such simple formulation is available for $\bm{Q}(\bm{z})$ in general. This problem was explored in \cite{Parsons2012}, where the following method was developed. This result is explained fully in Appendix \ref{Q}, for now we simply present the computational steps.

\begin{center}
\underline{Procedure for calculating $\bm{P}$ and $\bm{Q}$ at a point $\bm{z}\in\Gamma$}
\end{center}
\begin{enumerate}
\item Compute the Jacobian $\bm{J}$ and diagonalize it, writing
\begin{equation}
\bm{J}=\bm{W}\bm{\Lambda}\bm{W}^{-1}\,.
\label{diagJ}
\end{equation}
where $\bm{W}=(\bm{w}_1\cdots\bm{w}_n)$ is a matrix of eigenvectors forming a basis of $\mathbb{R}^d$, with the $m$ slow directions written first. $\bm{\Lambda}$ is a diagonal matrix of eigenvalues with $\lambda_1=\cdots=\lambda_m=0$ and $\textrm{Re}(\lambda_{m+1}),\ldots,\textrm{Re}(\lambda_n)<0$. Also compute the pseudo-inverse
\begin{equation}
\bm{J}^+=\bm{W}\bm{\Lambda}^+\bm{W}^{-1}\,,
\end{equation}
where $\bm{\Lambda}^+$ is the diagonal matrix with eigenvalues $\lambda^+_1=\cdots=\lambda^+_n$, where 
\[
	\lambda^+=\begin{cases}0\quad&\text{if}\quad \lambda=0\\1/\lambda&\text{if}\quad \lambda\neq 0\,.\end{cases}
\]
\item For each $i$, compute the Hessian $\bm{H}_i$ defined by $$H_{ijk}=\frac{\partial f_i(\bm{x})}{\partial x_j\partial x_k}\Bigg|_{\bm{x}=\bm{z}}\,.$$
Then find the (matrix-valued) solution $\bm{X}_i$ of the Lyapunov equation
\begin{equation}\label{lyapunov}
\bm{J}^T\bm{X}_i+\bm{X}_i\bm{J}=-\bm{H}_{i}\,.
\end{equation}
NB: this is a linear problem that is straightforwardly solved \cite{Bartels1972}.
\item Finally, the projection matrix is given by
\begin{equation}
\bm{P}=\bm{I}-\bm{J}^+\bm{J}\,.
\label{PmD}
\end{equation}
and for $\bm{Q}$ we have
\begin{equation}
Q_{ijk}=\sum_l -J^+_{il}[\bm{P^T}\bm{H}_l\bm{P}]_{jk}+P_{il}[\bm{X}_l-\bm{J}^{+T}\bm{H}_l\bm{P}-\bm{P}^T\bm{H}_l\bm{J}^+]_{jk}\,.
\label{QmD}
\end{equation}
\end{enumerate}

\subsection{Co-dimension one manifolds}

We now use the results of the previous section to obtain explicit expressions for the derivatives in the case when $\Gamma$ is a $(d-1)$-dimensional manifold. In this case, in a small neighbourhood around any point $\bm{z}\in\Gamma$, the flow field can be decomposed as $\bm{f} = \phi\,\bm{r}$, into a scalar part $\phi: \mathbb{R}^{d} \to \mathbb{R}$ that vanishes on $\Gamma$, and a non-vanishing vector part $\bm{r} : \mathbb{R}^{d} \to \mathbb{R}^{d}$. Using this decomposition we compute an expression for the Jacobian around a point:
\begin{equation*}
	\bm{J}(\bm{x}) = \phi(\bm{x}) \frac{\partial \bm{r}}{\partial \bm{x}}
		+ \bm{r}(\bm{x})\nabla \phi(\bm{x})^{T} .
\end{equation*}
In particular, evaluated exactly at the point $\bm{z}$ on the manifold we have $\bm{J} = \bm{r} \nabla\phi^{T}$. Meaning that $\bm{r}$ is, up to scalar multiple, the unique eigenvector corresponding to
\begin{equation*}
	\lambda = \nabla\phi^{T} \bm{r}\,,
\end{equation*}
which is the sole non-zero eigenvalue of $\bm{J}$. Note that here and hereafter we drop the arugment $\bm{z}$ to avoid notational clutter. As the Jacobian is given by the outer product of vectors $\bm{r}$ and $\nabla\phi$, it is straightforward to check that the pseudo-inverse may be written as
\begin{equation}
\bm{J}^+=\frac{1}{\lambda^2}\bm{J}^T\,.
\end{equation}
We conclude from \eqref{PmD} that
\begin{equation}\label{Pcd1}
\bm{P}=\bm{I}-\bm{J}^+\bm{J}=\bm{I}-\frac{\bm{J}^T\bm{J}}{\lambda^2}\,.
\end{equation}
To determine $\bm{Q}$, it thus remains to solve \eqref{lyapunov},
\begin{equation*}
\bm{J}^T\bm{X}_i+\bm{X}_i\bm{J}=-\bm{H}_{i}\,,
\end{equation*}
for $\bm{H}$ and insert into \eqref{QmD}. As developed in Appexdix \ref{Q}, equation \eqref{lyapunovsol}, the solution to this Lyapunov equation can be expressed by the exponential integral
\begin{equation*}
	\bm{X}_i
	= \int_{0}^{\infty} \Big(e^{s\bm{J}}-\bm{P}\Big)^{T} \frac{\partial^{2} f_{i}}{\partial \bm{x}^{2}} \Big(e^{s\bm{J}}-\bm{P}\Big)\,ds.
\end{equation*}
For an arbitrary vector $\bm{Y}$, we have $\bm{J} \bm{Y} = \bm{r} \nabla\phi^{T}\bm{Y}$, so that
\begin{equation*}
	\bm{J}^{2} \bm{Y}
	= \bm{r} \nabla\phi^{T}\bm{r} \nabla\phi^{T}\bm{Y}
	= \lambda (\nabla\phi^{T}\bm{Y}) \bm{r},
\end{equation*}
$\bm{J}^{n} \bm{Y} = \lambda^{n-1} (\nabla\phi^{T}\bm{Y}) \bm{r}$, and
\begin{equation*}
\begin{aligned}
	e^{s\bm{J}} \bm{Y} &= \sum_{n=0}^{\infty} \frac{s^{n}}{n!} \bm{J}^{n} \bm{Y}\\
	&= \bm{Y} + (\nabla\phi^{T}\bm{Y}) \bm{r} \lambda^{n-1}
		\sum_{n=1}^{\infty}  \frac{s^{n}}{n!}\\
	&= \frac{(\nabla\phi^{T}\bm{Y})}{\lambda}(e^{\lambda t}-1)
	\bm{r}.
\end{aligned}
\end{equation*}
Thus, $e^{s\bm{J}}-\bm{P} \bm{Y} =  \frac{(\nabla\phi^{T}\bm{Y})}{\lambda} e^{\lambda t} \bm{r}$ and, recalling that $\lambda < 0$,
we have that
\begin{equation*}
\begin{aligned}
	X_{ijk}
	&= \bm{e}_{j}^{T} \left(\int_{0}^{\infty} \bm{e}_{i}^{T} (e^{s\bm{J}}-\bm{P})^{T}
	\frac{\partial^{2} f_{i}}{\partial \bm{x}^{2}} (e^{s\bm{J}}-\bm{P})\,ds
	\right)\bm{e}_{k}\\
	&= \frac{(\nabla\phi^{T}\bm{e}_{j})(\nabla\phi^{T}\bm{e}_{k})}
	{\lambda^{2}}
	\bm{r}^{T} \frac{\partial^{2} f_{i}}{\partial \bm{x}^{2}} \bm{r}
	 \int_{0}^{\infty} e^{2 \lambda t}\, dt\\
	&= -\frac{\frac{\partial \phi}{\partial x_{j}}\frac{\partial \phi}{\partial x_{k}}}
	{2 \lambda^{3}}
	\bm{r}^{T} \frac{\partial^{2} f_{i}}{\partial \bm{x}^{2}} \bm{r}.
\end{aligned}
\end{equation*}
Finally, observing that $\frac{\partial^{2} f_{i}}{\partial x_{j}\partial x_{k}} =
\frac{\partial r_{i}}{\partial x_{j}}\frac{\partial \phi}{\partial x_{k}} + r_{i} \frac{\partial^{2} \phi}{\partial x_{j}\partial x_{k}}$, substituting the above into \eqref{QmD} and considerable algebraic simplification yields
\begin{multline}\label{Qcd1}
	Q_{ijk} = \frac{1}{\lambda} \left(
	[\bm{P}^{T}\frac{\partial^{2} \phi}{\partial \bm{x}^{2}}\bm{P}]_{jk}
		r_{i}
	- \frac{\partial \phi}{\partial x_{j}}
		[\bm{P} \frac{\partial \bm{r}}{\partial \bm{x}}]_{ik}
	- \frac{\partial \phi}{\partial x_{k}}
		 [\bm{P} \frac{\partial \bm{r}}{\partial \bm{x}}]_{ij}
	\right)+ \frac{1}{\lambda^{2}} \frac{\partial \phi}{\partial x_{j}}
		\frac{\partial \phi}{\partial x_{k}}
		[\bm{P} \frac{\partial \bm{r}}{\partial \bm{x}} \bm{r}]_{i}.
\end{multline}

\section{Worked examples}\label{WE}
\subsection{Simple example: Michaelis-Menten kinetics}
The Michaelis-Menten law is perhaps one of the most widely-applied examples of timescale separation. It is a model for the net rate of production in a chemical reaction that is catalysed by an enzyme, in which it is assumed that the process of enzyme binding and unbinding occurs very much faster than the catalytic reaction of interest. Using the notation of chemical reactions, one may write
\begin{equation}
E+S\xrightleftharpoons[k_r]{k_f} C \xrightarrow{k_{\text{cat}}} E+P\,,
\label{MM}
\end{equation}
where $E$ symbolises the enzyme, $S$ the substrate, $C$ the enzyme-substrate complex, and $P$ the product. The parameters $k_f$ and $k_r$ give the rate of binding (forward) and unbinding (reverse) of the enzyme to the substrate, while $k_{\text{cat}}$ specifies the rate of catalysis.

Assuming the reaction takes place in a domain of infinite volume, one may write rate the deterministic equations
\begin{equation}
\begin{split}
&\frac{dS}{dt}=-k_fES+k_rC\,,\\
&\frac{dE}{dt}=-k_fES+(k_r+k_{\text{cat}})C\,,\\
&\frac{dC}{dt}=k_fES-(k_r+k_{\text{cat}})C\,,\\
&\frac{dP}{dt}=k_{\text{cat}}C\,,\\
\end{split}
\end{equation}
where $S$, $C$, $P$ and $E$ now represent the \textit{concentrations} of the various reactants. Note that this system has only two degrees of freedom due to conservation relations $E+C=E_0$ and $S+C+P=S_0$, where $E_0$ and $S_0$ are the initial concentrations of the enzyme and substrate, respectively. If $k_f,k_r\gg k_{\text{cat}}$ we might approximate the concentration of the complex $C$ by the equilibrium value it would have if $k_{\text{cat}}$ were actually zero:
\begin{equation}
k_fES-k_rC\approx 0\quad \Rightarrow \quad C\approx E_0\frac{S}{k+S}\,,
\end{equation}
where $k=k_r/k_f$. Introducing $v^\ast=k_{\text{cat}}E_0$, on the slower timescale the net production rate is then found to be
\begin{equation}
\frac{dP}{dt}=\frac{v^\ast\, S}{k+S}\,.
\end{equation}
This is the Michaelis-Menten law.

In finite volume domains chemical reactions are subject to random fluctuations arising from the discrete nature of the molecules involved. A more appropriate description in these circumstances is a stochastic differential equation, with noise terms that are derived from the instantaneous reaction rates (each possible reaction introduces its own source of noise). For the reaction described above in (\ref{MM}) occurring in a domain of volume $V$, equations are derived following Kurtz \cite{Kurtz78}\footnote{In fact this step is not strictly necessary; we could choose to work directly with the process of particle numbers, as described in Appendix~\ref{Katz}.}:
\begin{equation}
\begin{split}
&\frac{dS}{dt}=-k_f(E_0-C)S+k_rC-\sqrt{\frac{k_f(E_0-C)S}{V}}\,\eta_{f}(t)+\sqrt{\frac{k_rC}{V}}\,\eta_{r}(t)\,,\\
&\frac{dC}{dt}=k_f(E_0-C)S-(k_r+k_{\text{cat}})C+\sqrt{\frac{k_f(E_0-C)S}{V}}\,\eta_{f}(t)-\sqrt{\frac{k_rC}{V}}\,\eta_{r}(t)-\sqrt{\frac{k_{\text{cat}}C}{V}}\,\eta_{\text{cat}}(t)\,.\\
%&\frac{dP}{dt}=k_{\text{cat}}C+\sqrt{\frac{k_{\text{cat}}C}{V}}\,\eta_{\text{cat}}(t)\\
\end{split}
\end{equation}
Following similar lines to \cite{Heineken1967} a dimensionless form may be found by rescaling time $t\mapsto k_fE_0t$ and introducing variables
\begin{equation}
\bm{x}=\left(\begin{array}{c}S/S_0\\C/E_0\end{array}\right)\,,
\end{equation}
and parameters
\begin{equation}
\begin{split}
&\varepsilon=\frac{k_{\text{cat}}}{k_fE_0}>0\,,\quad\mu=\frac{1}{S_0V}\,,\quad\alpha=\frac{k_r}{k_fS_0}>0 \,,\quad\beta=\frac{S_0}{E_0}>0\,.\\
\end{split}
\end{equation}
The result is a system of exactly the form of equation (\ref{sdex}):
\begin{equation}
\frac{d\bm{x}}{dt}=\bm{f}(\bm{x})+\varepsilon\bm{h}(\bm{x})+\sqrt{\mu}\,\bm{G}(\bm{x})\,\bm{\eta}(t)\,,
\end{equation}
where
\begin{equation}
\begin{split}
&\bm{f}(\bm{x})=\left(\begin{array}{c}-x_1+(x_1 +\alpha)x_2\\\beta(x_1-(x_1 +\alpha)x_2) \end{array}\right)\,,\qquad\bm{h}(\bm{x})=\left(\begin{array}{c}0\\-x_2 \end{array}\right)\,,\\
\end{split}
\label{MMfh}
\end{equation}
and
\begin{equation}
\begin{split}
&\bm{G}(\bm{x})=\left(\begin{array}{ccc}-\sqrt{(1-x_2)x_1}&\sqrt{\alpha x_2}&0\\\beta\sqrt{(1-x_2)x_1}&-\beta\sqrt{\alpha x_2}&-\sqrt{\varepsilon \beta x_2} \end{array}\right)\,,\qquad\bm{\eta}(t)=\left(\begin{array}{c}\eta_{f}(t)\\\eta_r(t)\\\eta_{\text{cat}}(t)\end{array}\right)\,.\\
\end{split}
\label{MMG}
\end{equation}
The slow manifold in this case is the curve $x_1-x_2(x_1+\alpha)=0$, along which $\bm{f}(\bm{x})=\bm{0}$. See Figure~\ref{MMfig} for an illustration.

Let us take $z=x_1$ as the slow variable and proceed to calculate a reduced system in terms of $z$ only. As the manifold is one-dimensional, we are able to simply follow the procedure laid out above. We begin by writing down the formula for the slow manifold and its $z$ derivatives:
\begin{equation}
\bm{\gamma}(z)=\left(\begin{array}{c}z\\\frac{z}{z+\alpha}\end{array}\right)\qquad \bm{\gamma}'(z)=\left(\begin{array}{c}1\\\frac{\alpha}{(z+\alpha)^2}\end{array}\right)\qquad\bm{\gamma}''(z)=\left(\begin{array}{c}0\\\frac{-2\alpha}{(z+\alpha)^3}\end{array}\right)\,.
\end{equation}
Next, we find the Jacobian matrix on the manifold
\begin{equation}
\bm{J}(\bm{x})=\left(\begin{array}{cc}x_2-1&x_1+\alpha\\\beta(1-x_2)&-\beta(x_1+\alpha) \end{array}\right)\quad\Rightarrow\quad \bm{J}(z)=\left(\begin{array}{cc}\frac{z}{z+\alpha}-1&z+\alpha\\\beta(1-\frac{z}{z+\alpha})&-\beta(z+\alpha) \end{array}\right)\,.
\end{equation}
Diagonalizing $\bm{J}(z)$ we find the left eigenvector $\bm{v}(z)$ corresponding to the eigenvalue zero, and its $z$ derivative:
\begin{equation}
\bm{v}(z)=\left(\begin{array}{c}\beta\\1\end{array}\right)\qquad \bm{v}'(z)=\left(\begin{array}{c}0\\0\end{array}\right)\,.\\
\end{equation}
Following equation (\ref{P1D}) we obtain
\begin{equation}
\bm{P}_{1\ast}(z)=\frac{(z+\alpha)^2}{\alpha+\beta(z+\alpha)^2}\left(\begin{array}{cc}\beta & 1\end{array}\right)\,,
\end{equation}
and from equation (\ref{Q1D})
\begin{equation}
\bm{Q}_{1\ast\ast}(z)=2\alpha\left(\frac{(z+\alpha)}{\alpha+\beta(z+\alpha)^2}\right)^3\left(\begin{array}{cc}\beta^2 & \beta\\ \beta&1\end{array}\right)\,.
\end{equation}
Plugging these results into the general formula (\ref{sdez2}) gives the reduced model
\begin{equation}
\frac{dz}{dt}=-\varepsilon\frac{z(z+\alpha)}{\alpha+\beta(z+\alpha)^2}+\varepsilon\mu\frac{\alpha\beta z(z+\alpha)^2 }{(\alpha+\beta(z+\alpha)^2)^3}-\frac{(z+\alpha)^2}{\alpha+\beta(z+\alpha)^2}\sqrt{\varepsilon\mu\frac{\beta z}{z+\alpha}}\,\eta_{\text{cat}}(t)\,.
\label{MMz}
\end{equation}
\begin{figure}
\includegraphics[width=0.6\textwidth, trim=0 0 0 0]{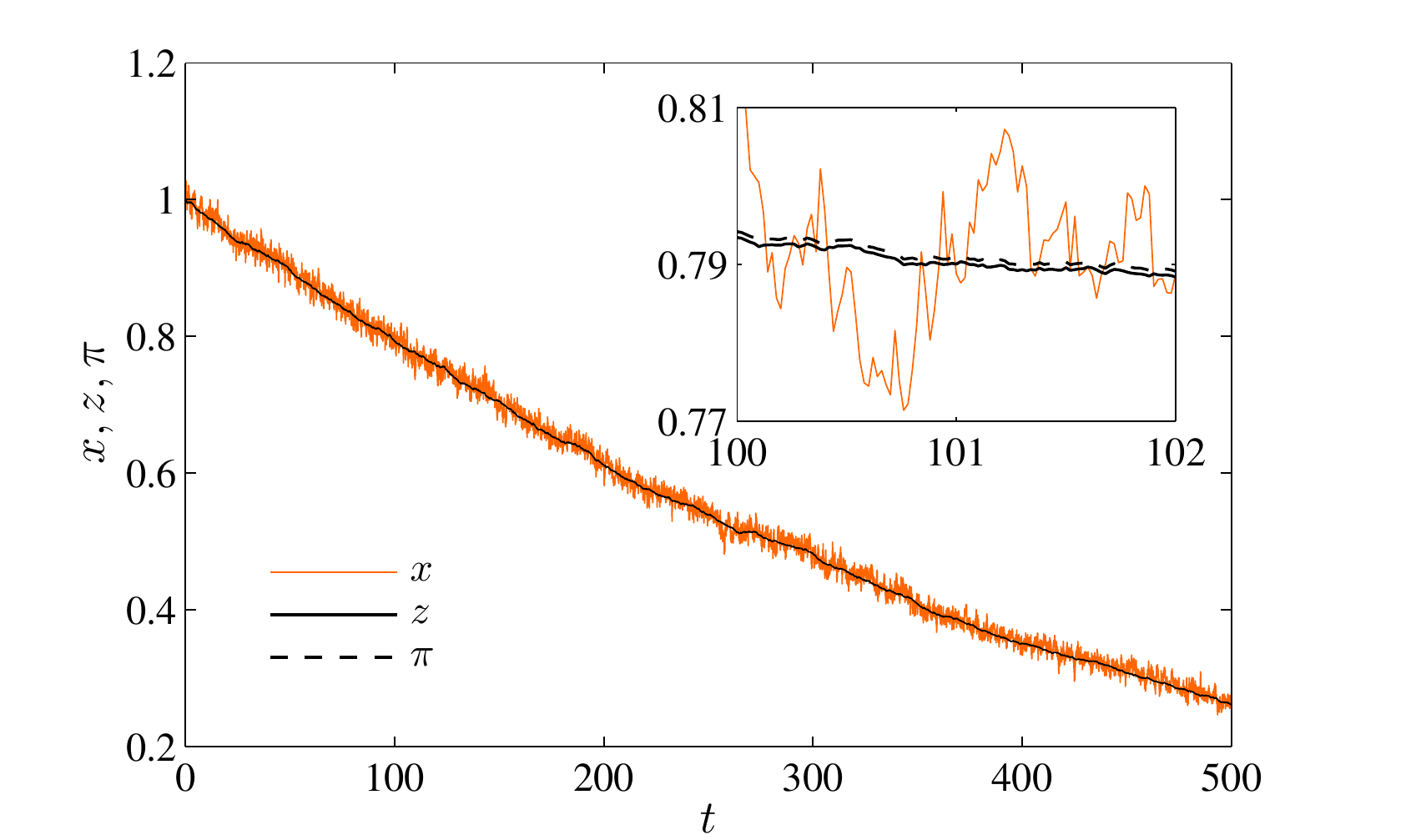}
\caption{Trajectories of $x_1$, $z$ and $\pi_1(\bm{x})$ from a single stochastic simulation of the Michaelis-Menten model. Note that the reduced dimension model for $z$ given by equation (\ref{MMz}) captures the dynamics of the full system under the projection $\bm{\pi}$ (hence the extremely close agreement between the solid and dashed black lines above). The original coordinate $x_1$ is subject to additional noise in the kernel of the projection. }
\label{MMfig2}
\end{figure}
Notice that there is a positive noise-induced drift term, meaning that the rate of decrease of $z$ is slowed by the noise. Figure~\ref{MMfig2} shows the dynamics of $z$ compared with those of $x_1$ in the full system for a single realization of the noise.

At first sight equation (\ref{MMz}) is considerably more complex than the traditional Michaelis-Menten law, however, carefully transforming back to the original coordinates\footnote{From the conservation rule $S+C+P=S_0$, we deduce that on the slow manifold we have $P=S_0-S_0z-E_0z/(z+\alpha)$. Use It\^o's lemma to compute $dP/dt$ from (\ref{MMz}) and finally undo the coordinate change via $t\mapsto t/k_fE_0$, $z\mapsto S/S_0$.}, one finds the simple result
\begin{equation}
\frac{dP}{dt}=\frac{v^\ast S}{k+S}+\sqrt{\frac{v^\ast S}{V(k+S)}}\,\eta_{\text{cat}}(t)\,.
\end{equation}

\subsection{Co-dimension one: the Wright-Fisher diffusion as a limit of a near-neutral stochastic Lotka-Volterra process}

Consider a well mixed-population of $d$ interacting species in an environment of carrying capacity $K$: there are $K$ ``slots'' in the environment that at most one individual may occupy. Let $X_{i}$ denote the number of individuals of species $i$, and suppose that each individual of species $i$ gives birth at rate $b_{i}$ and dies at rate $d_{i}$.  Further, suppose that the offspring is only viable if it lands in an empty patch, or if it lands in an occupied patch and out-competes the resident; say that an individual of type $i$ successfully displaces a resident of type $j$ with probability $c_{ij}$.  Then, there are three types of events:
\begin{enumerate}[(i)]
\item $X_{i}$ increases by 1 at rate $b_{i}X_{i} \left(1-\frac{\sum_{j} X_{j}}{K}\right)$,
\item $X_{i}$ decreases by 1 at rate $d_{i}X_{i}$, or,
\item $X_{i}$ increases by 1 and $X_{j}$ decreases by 1 at rate $b_{i} X_{i} \left(\frac{c_{ij}X_{j}}{K}\right)$.
\end{enumerate}
This gives a stochastic model of a population with density-dependent competition; \textit{n.b.} the total population size is \emph{not} fixed at $K$, but is rather allowed to fluctuate stochastically with an upper bound of $K$, as we allow the possibility of empty slots in the environment.

Let $x_{i}(t)$ denote the density of species $i$ (\textit{i.e.} $\frac{X_{i}(t)}{K}$).  As in the previous section, this system may be approximated by a system of stochastic differential equations, %Setting $b_{i} = b_{i}-d_{i}$ and $a_{ij} = b_{i} - b_{i} c_{ij}+b_{j}c_{ji}$ we obtain
\begin{multline*}
	\frac{dx_{i}}{dt} = \left((b_{i}-d_{i})- \sum_{j} (b_{i} - b_{i} c_{ij}+b_{j}c_{ji})x_{j}\right)x_{i}\\
	+ \sqrt{\frac{b_{i} x_{i}(1-\sum_{j} x_{j})}{K}} \eta_{b,i}(t)
	- \sqrt{\frac{d_{i} x_{i}}{K}} \eta_{d,i}(t)t
	+ \sum_{j} \sqrt{\frac{b_{i} c_{ij} x_{i} x_{j}}{K}} \eta_{i,j}(t)
	- \sum_{j} \sqrt{\frac{b_{j} c_{ji} x_{i} x_{j}}{K}} \eta_{j,i}(t).
\end{multline*}
%Thus, $\bm{f}(\bm{x})$ is of the Lotka-Volterra type discussed above.

We will be interested in large carrying capacity limits as $K \to \infty$, so here, $\mu = \frac{1}{K}$.
To explore the link between population genetics and population dynamics, we will further postulate that there exist values $\epsilon_{1},\ldots,\epsilon_{d}$ so that
\begin{equation*}
	b_{i} = b\left(1 + \frac{\epsilon_{i}}{K}\right), \quad
	d_{i} = d + \frac{\eta_{i}}{K},
	\quad \text{and} \quad c_{ij} = c + \frac{a_{ij}}{K},
\end{equation*}
for all $i,j$; this corresponds to the weak selection hypothesis \cite{Ewens1979}: all species are competitively equivalent and differ in their demographic rates by terms of $O\left(\frac{1}{K}\right)$. Then,
\begin{gather*}
	f_{i}(\bm{x}) =  x_{i}\left((b-d) - b \sum_{j=1}^{K} x_j\right),\\
	\varepsilon = \frac{1}{K},\\
\intertext{and}
	h_{i}(\bm{x}) = x_{i} \left(\left(b \epsilon_{i} - d \eta_{i}\right) -
		b \sum_{j} \left( (1-c)\epsilon_{i} - c \epsilon_{j} - a_{ij} + a_{ij} \right) x_{j} \right)
\end{gather*}
Under these assumptions,  $\text{rank}(\bm{A} \vert \bm{b}) = 1$,
\begin{equation*}
	\Gamma = \left\{\bm{x} \in \mathbb{R}^{d} : \sum_{j=1}^{K} x_{j} = 1-\frac{d}{b}\right\}
\end{equation*}
and for $\bm{x} \in \Gamma$, the derivatives \eqref{Pcd1} and \eqref{Qcd1} simplify to
\begin{equation*}
	P_{ij}(\bm{x}) = \delta_{ij} - \frac{x_{i}}{1-\frac{d}{b}} \quad \text{and} \quad
	Q_{ijk}(\bm{x})
		= - \frac{1}{1-\frac{d}{b}} \left(\delta_{ij} + \delta_{ik} - \frac{2 x_{i}}{1-\frac{d}{b}}\right),
\end{equation*}
whereas for $\bm{x} \in \Gamma$,
\[
	h_{i}(\bm{x}) = x_{i} \left(d(\epsilon_{i} - \eta_{i}) + c \sum_{j} (\epsilon_{i} - \epsilon_{j}) x_{j}
	  	+ b \sum_{j} (a_{ij}-a_{ji})x_{j}\right).
\]
A straightforward if lengthy calculation shows that $\bm{g}(\bm{x}) = O\left(\frac{1}{K^{2}}\right)$.

Substituting into our general formula \eqref{sdez2} then gives
\begin{multline*}
	\frac{d\tilde{x}_{i}}{dt}
	 = \frac{1}{K}
	 \left(h_{i}(\tilde{\bm{x}}) - \frac{\tilde{x}_{i}}{1-\frac{d}{b}} \sum_{j} h_{j}(\tilde{\bm{x}}) \right)\\
		+ \sum_{j} (\delta_{ij} - \frac{\tilde{x}_{i}}{1-\frac{d}{b}})
		\left(\sqrt{\frac{d \tilde{x}_{j}}{K}} (\eta_{b,j}(t) - \eta_{d,j}(t))
		+ \sum_{k} \sqrt{\frac{b c \tilde{x}_{j} \tilde{x}_{k}}{K}} (\eta_{j,k}(t) - \eta_{k,j}(t))\right),
\end{multline*}
or, changing variables to $p_{i} = \frac{\tilde{x}_{i}}{1-\frac{d}{b}}$, so $p_{i}$ is the proportion of species $i$,
\begin{multline*}
	\frac{dp_{i}}{dt}
	 = \frac{1}{K}
	 p_{i} \left(s_{i}(\bm{p}) -  \sum_{j} s_{j}(\bm{p}) p_{j} \right)\\
		+  \sum_{j} (\delta_{ij} - p_{i})
		\left(\frac{1}{\sqrt{1-\frac{d}{b}}} \sqrt{\frac{d p_{j}}{K}} (\eta_{b,j}(t) - \eta_{d,j}(t))
		+ \sum_{k} \sqrt{\frac{b c p_{j} p_{k}}{K}} (\eta_{j,k}(t) - \eta_{k,j}(t))\right),
\end{multline*}
where
\[
	s_{i}(\bm{p}) = d(\epsilon_{i} - \eta_{i})
		+ c\left(1 - \frac{d}{b}\right) \sum_{j} (\epsilon_{i} - \epsilon_{j}) p_{j}
	  	+ (b - d) \sum_{j} (a_{ij}-a_{ji})p_{j}.
\]

The corresponding Fokker-Planck equation for the density $f(\bm{p},t)$ is then
\begin{equation*}
		\frac{\partial f}{\partial t} = -\frac{1}{K} \frac{\partial}{\partial p_{i}} \left[
		p_{i} \left(s_{i}(\bm{p}) -  \sum_{j} s_{j}(\bm{p}) p_{j} \right) f\right]
		+ \frac{1}{2} \frac{{2\left(bc + \frac{d}{1-\frac{d}{b}}\right)}}{K}
		\frac{\partial^{2}}{\partial p_{i} \partial p_{j}}\left[p_{i}(\delta_{ij}-p_{j})f\right]
\end{equation*}
which we recognise as the equation for the Wright-Fisher diffusion, where the (frequency dependent) selection coefficient is $\frac{s_{i}(\bm{p})}{K}$ and the effective population size is $N_{e} = \frac{\left(1-\frac{d}{b}\right)K}{2\left(c(b-d) + d\right)}$; $\left(1-\frac{d}{b}\right)K$ is the population size at the deterministic equilibrium, whereas the other terms reflect variance in the total population size.  %Analogous expressions are obtained when one considers $O\left(\frac{1}{K}\right)$ differences in the values $d_{i}$ or $c_{ij}$.
This gives an alternate derivation of the results presented in \cite{Parsons2007b,Parsons2008,Parsons2010,Constable2015}.

\subsection{Continuous degrees of freedom: example of competition-limited diffusion}
The methods of Section \ref{S2} can readily be extended to infinite dimensional settings. Two recent examples come from work exploring the role of stochasticity in spatial ecological models \cite{Rogers2012a,Pigolotti2014}. Here we work through a simple illustrative example of diffusing particles coupled by a competitive birth-death interaction; we will show that this competition acts to limit the speed of diffusion of the population. Interested readers are referred to \cite{Etheridge1991}, where the continuum limit of this example has been studied in considerable depth.

Consider the following stochastic process. At time $t$ there are $N_t$ individual particles wandering in a one-dimensional space, each following their own Brownian motion with diffusion constant $D=\sqrt{2\varepsilon}$. With rate one, each particle may independently `reproduce', creating a daughter particle that initially shares the location of the parent, but thereafter moves independently. Particles `die' with rate proportional to their total number; specifically, the death rate for each particle is $\mu(N_t-1)$. We assume the constants $\mu$ and $\varepsilon$ are small, but of the same order.

\begin{figure}
\includegraphics[width=0.8\textwidth, trim=50 0 50 0]{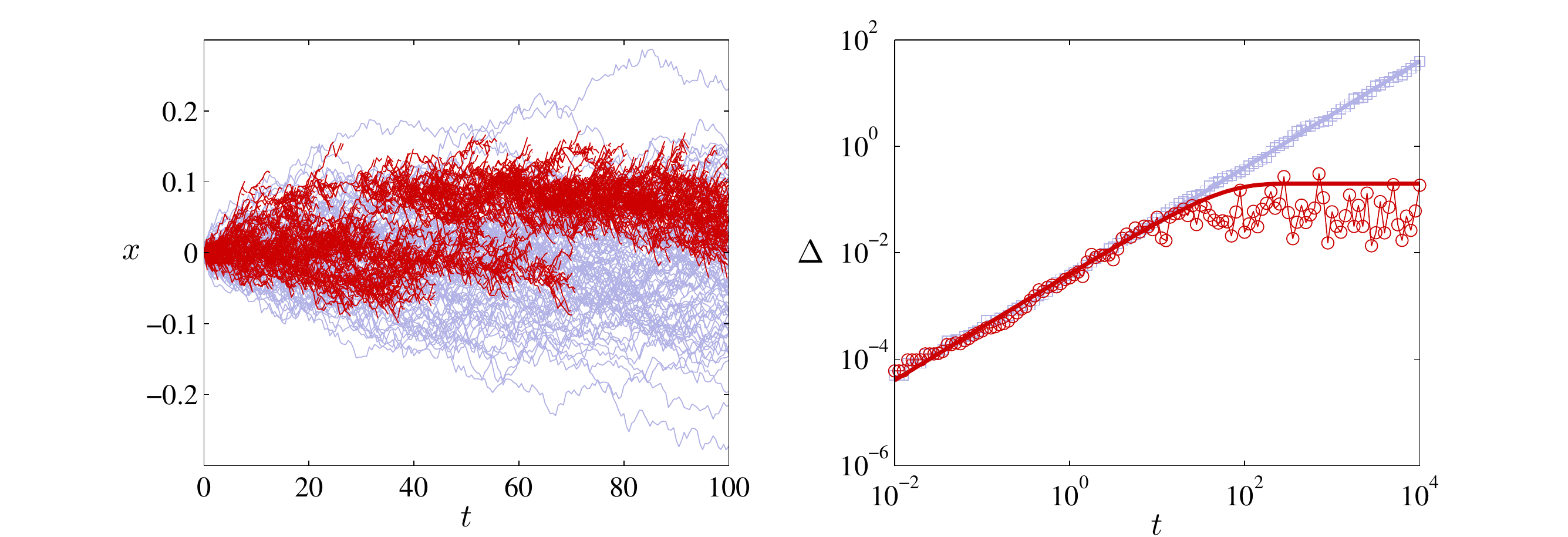}
\caption{Simulation of competition-limited diffusion (dark red), contrasted with a collection of $N$ independent Brownian particles (light purple). The left panel shows the particle trajectories, on the right is shown the mean square distance between pairs of particles. }
\label{clbm_fig}
\end{figure}

Since the location of the particles does not influence the birth or death rates, it is easy to see that the total number of particles follows a logistic growth law, quickly reaching an equilibrium $N_t\approx \mu^{-1}$. The total population size remains at this level while the spatial distribution of particles evolves slowly over a much longer timescale. We are interested in the long-term behaviour of the distribution of particle locations. Introduce the population density
\begin{equation}
u(x)=\mu\sum_{n=1}^{N_t}\delta\left(x-X_n(t)\right)\,,
\end{equation}
where $X_n$ is the location of particle $n$ at time $t$, $\delta$ is the Dirac delta function, and we suppress the dependence of $u$ on $t$ to reduce clutter. Simulations suggest that the competitive interaction of the particles limits the extent to which they are able to diffuse away from each other (Figure~\ref{clbm_fig}, left panel). This observation can be made quantitative by computing the mean square distance between pairs of particles,
\begin{equation}
\Delta[u]:=\mu^2\sum_{n,m}(X_n(t)-X_m(t))^2=\iint (x-y)^2u(x)u(y)\,dx\,dy\,.
\end{equation}
The right panel of Figure~\ref{clbm_fig} shows the time evolution of $\Delta$ for the population, compared to the growth $\Delta \sim t$ observed for independent diffusing particles. The solid lines show our theoretical prediction for this phenomenon, which we will now derive using timescale separation.

Following a system-size expansion \cite{McKane2014}, we find that the time-evolution of $u(x)$ is described to close approximation by the stochastic partial differential equation (SPDE)
\begin{equation}
\frac{\partial}{\partial t}u(x)=\varepsilon\frac{\partial^2}{\partial x^2}u(x)+u(x)\left(1-\int u(y)\,dy\right)+\sqrt{\mu\, u(x)\left(1+\int u(y)\,dy\right)}\,\eta(x,t)\,,
\label{spdeu}
\end{equation}
where $\eta(x,t)$ is spatio-temporal white noise and the integrals run over the real line.

Equation (\ref{spdeu}) has the same essential structure as our basic object of interest (\ref{sdex}). If we identify
\begin{equation}
\begin{split}
&f[u](x)=u(x)\left(1-\int u(y)\,dy\right)\\
&h[u](x)=\frac{\partial^2}{\partial x^2}u(x)\\
&G[u](x,s)=\delta(x-s)\sqrt{u(x)\left(1+\int u(y)\,dy\right)}\,,
\end{split}
\label{fhG_ex}
\end{equation}
then (\ref{spde}) becomes
\begin{equation}
\frac{\partial }{\partial t}u(x)=f[u](x)+\varepsilon h[u](x)+\sqrt{\mu}\int G[u](x,s)\eta(s,t)\,ds\,.
\label{spde}
\end{equation}
The integral here is the analogue of the matrix-vector multiplication $\bm{G}(\bm{x})\bm{\eta}(t)$ appearing in (\ref{sdex}). The delta function appearing in $G[u]$ means that the noise in our example is spatially uncorrelated; this may not hold for other models.

In this section we will show how the timescale separation techniques discussed above may also be applied to equations of the form (\ref{spde}). First, an important caveat: the irregular nature of spatio-temporal noise creates enormous mathematical complications in the rigorous analysis of SPDEs, we refer interested readers to the Fields-medal winning work \cite{Hairer2014}. In what follows we will turn a blind eye to deeper questions concerning the nature of the solution space and simply apply the techniques developed in the previous sections.

First we examine the outer part $\partial u/\partial t=f[u]$. In our example, the PDE
\begin{equation}
\frac{\partial}{\partial t}u(x)=u(x)\left(1-\int u(y)\,dy\right)\,,
\end{equation}
is straightforward to solve:
\begin{equation}
u(x,t)=\frac{u(x,0)\,e^t}{1+(e^t-1)\int u(y,0)\,dy}\,,
\end{equation}
which describes the fast relaxation of $u$ to a state in which it has total mass one. In this infinite-dimensional setting, the map that describes the long-time limit of the outer solution (previously defined in (\ref{defpi})) is an operator $\pi$, whose action is specified by
\begin{equation}
\pi[u](x)=\frac{u(x)}{\int u(y)\,dy}\,.
\end{equation}
We suppose that there exists a suitable space of functions $\mathcal{U}$ describing possible solutions of (\ref{spde}). Exactly what kind of space is a deep question beyond our present focus. The analogue of the slow manifold is the subspace $\mathcal{V}\subset\mathcal{U}$ containing functions $v$ satisfying $f[v]=0$, or equivalently for our example, $\int v(y)\,dy = 1$. We aim to derive an equation describing slow stochastic evolution in $\mathcal{V}$ that well-approximates the behaviour of solutions to the full system (\ref{spde}).

Where previous calculations involved partial differentiation, we now apply a functional derivative. In analogue to the definitions in (\ref{defPQ}) we introduce
\begin{equation}
P[v](x,y)=\frac{\delta }{\delta u(y)}\pi[u](x)\Bigg|_{u=v}\,,\quad Q[v](x,y,z)=\frac{\delta^2 }{\delta u(y)\delta u(z)}\pi[u](x)\Bigg|_{u=v}\,.
\end{equation}
The reduced system may then be written down:
\begin{equation}
\begin{split}
\frac{\partial}{\partial t}v(x)=&\int P[v](x,y)\left[\varepsilon h(y)\,dy +\sqrt{\mu}\int G[v](y,s)\eta(s,t)\,ds\right]\\
&+\frac{\mu}{2}\iiint G[v](y,s)G[v](z,s)Q[v](x,y,z)\,dy\,dz\,ds\,.
\end{split}
\label{spde_re}
\end{equation}

For the example at hand we compute
\begin{equation}
\begin{split}
&\frac{\delta }{\delta u(y)}\pi[u](x)=\frac{\delta(x-y)}{\int u(z)\,dz}-\frac{u(x)}{\left(\int u(z)\,dz\right)^2}\,,\\
&\frac{\delta^2 }{\delta u(y)^2}\pi[u](x)=\frac{2u(x)}{\left(\int u(z)\,dz\right)^3}-\frac{2\delta(x-y)}{\left(\int u(z)\,dz\right)^2}\,,
\end{split}
\end{equation}
and thus
\begin{equation}
P[v](x,y)=\delta(x-y)-v(x)\,,\quad Q[v](x,y,y)=2v(x)-2\delta(x-y)\,.
\label{PQ_ex}
\end{equation}
Note that we only need the $z=y$ parts of $Q[v](x,y,z)$ because of the delta function in $G$. Plugging (\ref{fhG_ex}) and (\ref{PQ_ex}) into (\ref{spde_re}), we obtain the reduced model
\begin{equation}
\frac{\partial}{\partial t}v(x)=\varepsilon \frac{\partial^2}{\partial x^2}v(x)+\sqrt{2\mu}\int \big[\delta(x-y)-v(x)\big]\sqrt{v(y)}\eta(y,t)\,dy\,.
\label{spdev}
\end{equation}

Comparing (\ref{spdev}) to the original equation (\ref{spdeu}) we see two main differences: the non-linearity in the drift has vanished, but the noise is now spatially coupled.

To compute a prediction for the mean squared distance between particles, it is simpler to work in Fourier space. Introducing $\tilde{v}(k)=\int e^{-2\pi i k x}v(x)\,dx$, we note first that
\begin{equation}
\mathbb{E}\Delta[v]=\iint z^2e^{2\pi i k z}\,\mathbb{E}|\tilde{v}(k)|^2\,dk\,dz=-\frac{1}{4\pi^2}\frac{\partial^2}{\partial k^2}\mathbb{E}|\tilde{v}(k)|^2\Bigg|_{k=0}
\label{EDv}
\end{equation}
Translating (\ref{spdev}) to Fourier space we find
\begin{equation}
\frac{\partial}{\partial t}\tilde{v}(k)=-4\varepsilon\pi^2k^2\tilde{v}(k)+\sqrt{2\mu}\int \widetilde{G}[\tilde{v}](k,x)\,\eta(x,t)\,dx\,,
\end{equation}
where
\begin{equation}
\widetilde{G}[\tilde{v}](k,x)=\left(e^{-2\pi ikx}-\tilde{v}(k)\right)\sqrt{\int\,e^{2\pi i \ell x}\tilde{v}(\ell)d\ell}\,.
\end{equation}
In mean, this process behaves exactly as a straightforward diffusion:
\begin{equation}
\frac{d}{dt}\mathbb{E}[\tilde{v}(k)]=-4\varepsilon\pi^2k^2\,\mathbb{E}[\tilde{v}(k)]\,.
\end{equation}
However, the noise introduces a correction to the variance following It\^o's formula. Specifically,
\begin{equation}
\begin{split}
\frac{d}{dt}\mathbb{E}|\tilde{v}(k)|^2&=-8\varepsilon\pi^2k^2\mathbb{E}|\tilde{v}(k)|^2+\frac{1}{2}\iiint \widetilde{G}[\tilde{v}](\ell,x)\widetilde{G}[\tilde{v}](m,x)\frac{\delta^2|\tilde{v}(k)|^2}{\delta \tilde{v}(\ell)\delta \tilde{v}(m)} \,dx\,d\ell\,dm\\
&=-8\varepsilon\pi^2k^2\mathbb{E}|\tilde{v}(k)|^2+2\mu(1-\mathbb{E}|\tilde{v}(k)|^2)\,.
\label{dEv2}
\end{split}
\end{equation}
Solving (\ref{dEv2}) and plugging into (\ref{EDv}) gives the prediction
\begin{equation}
\mathbb{E}\Delta[v]=\frac{2\varepsilon}{\mu}\left(1-e^{-2\mu t}\right)\,.
\end{equation}
This result is shown as the dark red curve in Figure ~\ref{clbm_fig}. In particular, notice that  whilst the mean square distance between diffusing particles grows indefinitely, in the competition coupled process it attains a finite limit $2\varepsilon/\mu$.

\section{Discussion}

The purpose of this article has been to show the derivation and application of a systematic computational framework for dimension reduction in stochastic dynamical systems that exhibit a separation of timescales via a globally stable normal hyperbolic slow manifold \textit{i.e.}\, in the limit of small noise the limiting deterministic dynamical system defined by $\bm{f}$ possesses a single, connected and globally attractive manifold of fixed points.   The method is exact in the limit of small noise and well-separated slow and fast dynamics, and experimentally found to be valid as an approximation scheme over a sensible parameter range. We have also presented extensions of the method for infinite dimensional systems and processes coupled to general noise sources. 

%Our analysis rests on some assumptions about the behaviour of the outer system; we have focused here solely on the situation that the deterministic dynamical system defined by $\bm{f}$ possesses a single, connected and globally attractive manifold of fixed points. 

In some applications more general scenarios may occur, we now briefly discuss two of interest. Some models may exhibit more than one connected manifold of equilibria or dynamic bifurcations, \text{i.e.}, points where the critical manifold ceases to be normally hyperbolic \cite{Berglund2006}; in this case the theory developed here will apply locally to trajectories in the basin of attraction of each manifold individually, but further analysis will be necessary to describe the statistics of noise-driven transitions between manifolds. A possibly more exciting direction for further research is the analysis of noisy behaviour around more general attractors such as limit cycles, limit tori and strange attractors. In the case of limit cycles some work exists on stochastic extensions to Floquet theory \cite{Boland2009}, however, this is a linear description that cannot capture any bias analogous to the noise-induced drift in the slow manifold setting. 

Finally, it is worth returning to discuss the motivation for this work. As mentioned earlier, variations of the work of Katzenberger have been independently rediscovered by several groups in recent years, almost all of whom have been interested in questions about the role of noise in ecology and evolution. Historically, many theoretical results in this field have been derived from models that assume for convenience a fixed population size. In the deterministic limit this assumption is not important, but we are now beginning to realise that the inclusion of noise can induce radically different and sometimes unexpected behaviour. Mathematically, this is a consequence of the noise-induced drift term $\bm{g}$ that appears in our equation (\ref{sdez}), and more generally of the seemingly endless capacity of It\^{o}'s lemma to cause surprise. There have been some tentative explorations of the possible evolutionary and ecological consequences of these effects \cite{Parsons2010,Rossberg2013,Constable2014b}, but much more is yet to be discovered.

\section*{Acknowledgements} We thank George Constable, Alan McKane and Christopher Quince for useful discussions and Tom Kurtz for bringing \cite{Katzenberger91} to our attention.  Work on this paper began when both TLP and TR were visitors at the Isaac Newton Institute for Mathematical Sciences programme \textit{Understanding Microbial Communities; Function, Structure and Dynamics}. 

\section*{Funding statement} TLP is supported by the CNRS, and TR by the Royal Society. 

\appendix
\renewcommand{\thesubsection}{\thesection.\arabic{subsection}}

\section{Related literature}\label{lit}
As mentioned in the introduction, the mathematical and theoretical physics literatures contain a multitude of techniques for separation of timescales, many of which have been extremely well-studied. In this appendix we present a brief overview of some of the historical developments that we consider to be more relevant to the class of systems we are interested in.

The canonical example of timescale separation in physics is Brownian motion; stochastic interactions with the water molecules cause changes to the \emph{velocity} of the pollen grain on a fast timescale, which have a cumulative effect of perturbing the \emph{position} on a slower timescale. The process of moving from a description of the particle's motion in terms of position and velocity to one concerned just with position (thus reducing the dimension of the model from two to one) is known as \textit{adiabatic elimination}\footnote{Also \textit{fast variable elimination}, \textit{fast mode elimination}, the \textit{quasi-static approximation}, and many other alternatives.}. Following the introductory discussion in \cite{Gardiner1985}, one may write equations of motion for the pollen grain of the form of Langevin equations such as
\begin{equation*}
\frac{dx}{dt}=v\,,\quad m\frac{dv}{dt}=-v+\sqrt{\mu}\,\eta(t)\,,
\end{equation*}
where $m$ is the mass, $\mu$ a constant derived from the temperature of the water bath, and $\eta(t)$ is Gaussian white noise. If the mass $m$ is very small, we might approximate $m\approx0$ and thus $v\approx \sqrt{\mu}\,\eta(t)$. From this we derive the reduced description
\begin{equation*}
\frac{dx}{dt}=\sqrt{\mu}\,\eta(t)\,.
\end{equation*}
A reader encountering arguments of this type for the first time is likely to be suspicious, and rightly so. Although provably exact for the simple case described here, taking limits in a brusque fashion like this is generally inadvisable when dealing with stochastic dynamical systems. Nonetheless, the method is quite powerful and the same basic principle has been developed to various levels of rigour and generality by many authors, most notably Haken \cite{Haken1982}.

Although conceptually appealing, the Langevin-style description of a stochastic system in terms of dynamical equations with noise hides some important complications. As any introductory text will warn, when writing such an equation we must specify the sense in which we are to understand the noise term $\eta(t)$. This choice then impacts the behaviour of the system under certain limits or changes of variables, potentially complicating the process of timescale separation. An alternative description of a stochastic system that avoids this ambiguity is in terms of a PDE describing the time evolution of the probability density, known as the Fokker-Planck Equation (also, the Kolmogorov forward equation). Timescale separation in this setting amounts to integrating out one or more degrees of freedom from the PDE to reduce its dimension. Physicists might understand this process through its natural analogue in quantum mechanics, the Born-Oppenheimer approximation \cite{Born1927}. In most applications the Fokker-Planck equation will not be exactly separable, necessitating the application by hand of a carefully chosen projection operator, an approach going back to the work of Zwanzig \cite{Zwanzig1960}.

Existing in parallel with the development of adiabatic elimination techniques in physics is a separate body of literature concerned with probabilistic models of biological processes that also exhibit a separation of timescales. This thread begins with early work on the convergence of Markov chains --- particularly those appearing in genetics --- to diffusion processes \cite{Feller1951,Trotter1958}. Analyzing the Wright-Fisher model of allele frequencies in a population of large size $N$, Feller  used the standard convention of measuring the population in units of $N$ individuals (so one individual has weight $\frac{1}{N}$), but rather than truncating the master equation by discarding terms of order $\frac{1}{N^{2}}$ and greater, he also rescaled to a ``slow time'' so that $\Delta t = \frac{1}{N}$, which led to a well defined limit as $N \to \infty$, a procedure proved rigorously by Trotter.

Mathematical population genetics also provides an important early example of a stochastic process explicitly considered at multiple timescales. In \cite{Ethier1980,Nagylaki1980}, models of diploid populations with non-overlapping generations and geographic structure were analysed by separating a fast timescale on which, for example, the genotype frequencies would rapidly equilibrate to Hardy-Weinberg proportions, from a slow timescale over which allele frequencies would vary due to mutation and genetic drift. Their approach, obtained by an extension of the semi-group methods of \cite{Trotter1958}, requires the explicit characterization of infinite-dimensional spectral projection operators, and an explicit separation of the process into fast and slow variables, the former of which must be effectively constant on the slow-timescale.  These requirements somewhat limit the generality of the approach, but it does have the virtue of being applicable to Markov processes other than Brownian motion, such as measure-valued processes \cite{Harris2015}.

Slightly over a decade later, two very general papers appeared that used  similar approaches based in stochastic differential equations and classical, deterministic formulations of time-scale separation.  The first, \cite{Katzenberger91}, considered processes taking values in $\mathbb{R}^{d}$, but with a very general approach to stochastic noise, so that the results can be applied both to diffusion processes and discrete Markov processes (see the discussion in section \ref{Katz}), whilst the latter \cite{Funaki1993}, allows the process to be defined on a general Riemann manifold $M$, but requires the noise to be the canonical Brownian motion on the manifold. %(see section \ref{Funaki}). 
In both, there is a fast timescale on which the process essentially behaves like a deterministic dynamical system, with trajectories that approach a lower-dimensional ``slow manifold'' $\Gamma$.  On the slow time-scale, the process is asymptotic to a diffusion process that is confined to the manifold $\Gamma$. Katzenberger characterized the diffusion on the slow manifold via stochastic differential equations and a function $\bm{\pi}$ (defined formally below), which, given an initial point $\bm{x}$, gives the point in $\Gamma$ to which it will be carried by the fast dynamics, whereas \cite{Funaki1993} used the backward equation to describe the diffusion.

Unfortunately, in spite of their generality, \cite{Katzenberger91} and \cite{Funaki1993} have been frequently overlooked, as was the follow up \cite{Funaki1995}, where the problem of dimension reduction in infinite dimensions was considered for an SPDE of Ginzburg-Landau type. Indeed, special cases have been subsequently rediscovered.  This is perhaps because, whilst these papers have the virtue of providing rigorous proofs, they are not necessarily useful for the practitioner:  \cite{Katzenberger91} supposes the function $\bm{\pi}$ is given, whereas in applications it is frequently difficult or impossible to obtain in closed form, as the fast system may not be solvable.  \cite{Funaki1993} gives a closed expression for the diffusion equation, though one that, except in the case when the fast system is a gradient flow, requires extensive calculations with Fermi coordinates (like Katzenberger's $\bm{\pi}$, these are often impossible to obtain in closed form) and only applies when the noise is Brownian motion.

\section{Katzenberger's Theorem}\label{Katz}

Above we developed our results in the context of It\^o SDEs, however, \cite{Katzenberger91} proved a more general result that allows us to consider a much broader class of noise processes: semimartingales. Semimartingales are the most general class of stochastic processes for which one may define a stochastic integral and stochastic differential equations (Brownian motion is included as a special case). Suitably adapted, most of the familiar results for SDEs and white-noise integrals, including It\^o's formula, remain true in the more general setting \cite{Protter04}.

To define a semimartingale, we must first make a few auxiliary definitions.  A Markov process $M(t)$ is a \textit{martingale} if
\begin{equation*}
	\mathbb{E}\left[M(t) \middle\vert M(s) \right] = M(s).
\end{equation*}
A random variable $\tau$ taking values in $[0,\infty)$ is a \textit{stopping time} if one can determine if $\tau < t$ without knowledge of the future beyond $t$; an example of a stopping time is the first time a diffusion started from 0 exits an interval $[-a,a]$. $M(t)$ is a \textit{local martingale} if there is a sequence of stopping times $\tau_{n} \to \infty$ such that $M(\min\{t,\tau_{n}\})$ is a martingale for each $n$.

A function is \textit{c\`adl\`ag} if it is continuous from the right and has left-hand limits at every point.

The \textit{total variation} of a function $f$ on an interval $[a,b]$ is
\[
	V^{a}_{b}(f) = \min_{\{t_{i}\}} \sum_{i} |f(t_{i+1})-f(t_{i})|,
\]
where the minimum is over all partitions $a=t_{0} < t_{1} < \cdots < t_{n} = b$ of $[a,b]$. A stochastic process $A(t)$ is of \textit{finite variation} if it is c\`adl\`ag and has finite total variation on all intervals $[a,b]$ (note that $A(t)$ is allowed to have jump discontinuities).

Finally, $Z(t)$ is a semimartingale if it may be written as the sum of a local martingale and a finite variation process,
\[
	Z(t) = M(t) + A(t).
\]
Diffusion processes are the prototypical example of semimartingales, but the class is much broader, and includes processes with jumps, such as L\'evy processes; \textit{e.g.}\, if $N(t)$ is a Poisson process, then $M(t) = N(t) - t$ is a local martingale and $A(t) = t$ is of finite variation, so $N(t)$ is a semimartingale.  Integration with respect to a semimartingale is defined analogously to the Stieltjes integral, except that we require the approximating sum to converge in probability, and, as with the It\^o integral, the integrand is always evaluated at the left endpoint of each interval in the partition.

More generally, we can define vector and matrix valued martingales, local martingales, finite variation processes and semimartingales, $\bm{M}(t)$, $\bm{A}(t)$, and $\bm{Z}(t)$, by requiring the components, $M_{i}(t)$ \textit{etc.}, have the corresponding property.

We can now formulate Katzenberger's result.  Let
\begin{enumerate}
\item[(i)] $\bm{Z}_{n}(t)$ be a convergent sequence of vector valued semimartingales such that the jumps $\Delta \bm{Z}_{n}(t) \to 0$ as $n \to \infty$,
\item[(ii)] $A_{n}(t)$ be a sequence of non-decreasing finite variation processes such that $\Delta A_{n}(t) \to 0$, and
\[
	\int_{a}^{b} dA_{n}(s) = A_{n}(b) - A_{n}(a) \to \infty;
\]
Katzenberger notes that most frequently in applications, $A_{n}(t) = \alpha_{n} t$ for some sequence $\alpha_{n} \to \infty$ (\textit{n.b.}, in this formulation, this explosion in $A_{n}(s)$ corresponds to the drift becoming infinitely strong, rather than the noise infinitely weak, as in \eqref{sdex}.   The two are equivalent, if one changes the timescale accordingly; recall we had
\begin{equation*}
\frac{d\bm{x}}{dt}=\bm{f}(\bm{x})+\varepsilon\bm{h}(\bm{x})+\sqrt{\mu}\,\bm{G}(\bm{x})\,\bm{\eta}(t).
\end{equation*}
If instead, we consider the process $\tilde{\bm{x}}(t) = \bm{x}(\mu t)$, we get
\begin{equation*}
\frac{d\tilde{\bm{x}}}{dt}= \frac{1}{\mu} \bm{f}(\tilde{\bm{x}})+\frac{\varepsilon}{\mu}\bm{h}(\tilde{\bm{x}})+\,\bm{G}(\tilde{\bm{x}})\,\bm{\eta}(t),
\end{equation*}
with a drift that blows up as $\mu \to 0$).
\item[(iii)] $\bm{f}$ and $\Gamma$ be as before,
\item[(iv)]  $\bm{G}_{n}(\bm{x})$ be a sequence of matrix-valued functions converging to a limit $\bm{G}(\bm{x})$, and
\item[(v)] $\bm{x}_{n}(t)$ be a sequence of stochastic processes satisfying the (semimartingale) SDE
\[
\	d\bm{x}_{n} = \bm{f}(\bm{x}_{n})\, dA_{n} + \bm{G}_{n}(\bm{x}_{n})\, d\bm{Z}_{n}.
\]
\end{enumerate}
Then, as before, subject to a few technical considerations, as $n \to \infty$, $\bm{x}_{n}$ converges to a diffusion process on $\Gamma$ satisfying %the It\^o SDE
%\[
%	d\bm{z} = \bm{g}(\bm{z})\, dt + \bm{P}(\bm{z})\bm{G}(\bm{z})\, d\bm{B}(t),
%\]
%or, equivalently, in Langevin notation,
\begin{equation}\label{KSDE}
\frac{d\bm{z}}{dt}=\bm{g}(\bm{z})+ \bm{P}(\bm{z})\bm{G}(\bm{z})\bm{\eta}(t)\,,
\end{equation}
where $\bm{g}$ is as in equation \eqref{drift} and $\bm{\eta}$ is white noise.

\begin{figure}
\includegraphics[width=0.8\textwidth, trim=50 320 50 280]{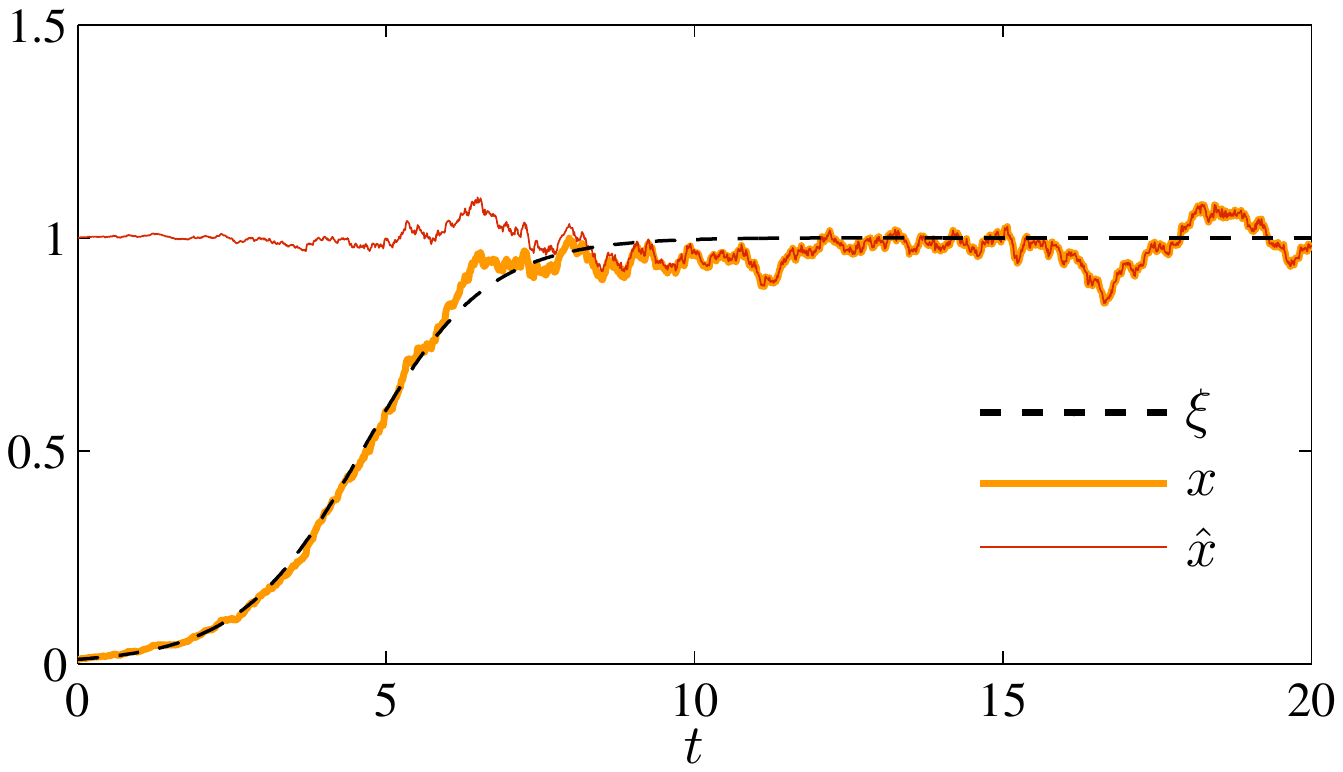}
\caption{Illustration of $\hat{x_1}$ on the fast timescale for a prototypical stochastic dynamical system with a slow manifold $\Gamma=\{\bm{x}\,:\,x_1=1\}$. In the slow timescale the initial transit to the manifold is compressed into an instantaneous jump at $t=0$.}
\label{xi_fig}
\end{figure}

Some care is required in understanding the sense of convergence in \cite{Katzenberger91}; if
$\bm{x}_{n}(0)$ converges weakly to $\bm{z} \in \Gamma$ in $\mathbb{R}^{d}$ (\textit{i.e.},\, for all continuous functions $F:\mathbb{R}^{d} \to \mathbb{R}$, $\mathbb{E}[F(\bm{x}_{n}(0))] \to \mathbb{E}[F(\bm{z})]$)
then $\bm{x}_{n}(t)$ converges weakly to $\bm{z}(t)$ in the space of c\`adl\`ag functions:
\[
	\mathbb{E}[F(\bm{x}_{n}(t))] \to \mathbb{E}[F(\bm{z}(t))]
\]
for all continuous functions $F$ from the space of c\`adl\`ag functions on $[0,\infty)$ to $\mathbb{R}$, see \cite{Billingsley68,Ethier+Kurtz86} for a definition of the topology on c\`adl\`ag functions and results on weak convergence.  When $\bm{x}_{n}(0)$ converges to a limit  $\bm{x}$ that is not in $\Gamma$, additional care is required: in this case, the process will jump instantaneously from $\bm{x}$ to $\bm{\pi}(\bm{x}) \in \Gamma$, which is not compatible with convergence in the weak topology on c\`adl\`ag functions.  However, if one considers
\[
	\hat{\bm{x}}_{n}(t) = \bm{x}_{n}(t) - \bm{\xi}(nt) + \bm{\pi}(\bm{x}),
\]
(recall, $\bm{\xi}(t)$ is the solution to the outer system, \eqref{ivp}) then
$\hat{\bm{x}}_{n}(0) \to \bm{\pi}(\bm{x}) \in \Gamma$ and $\hat{\bm{x}}_{n}(t)$ converges weakly to the diffusion $\bm{z}(t)$ on $\Gamma$  as before; intuitively $\hat{\bm{x}}_{n}(t)$ is obtained by removing the initial transient phase when $\bm{x}_{n}(t)$ follows the trajectories of the outer system, and starting the process instead from the endpoint of that trajectory, $\bm{\pi}(\bm{x})$ (see Figure \ref{xi_fig}).

\subsection{Density dependent population processes}

While Katzenberger's result might seem unnecessarily abstract, it allows one to apply the same slow-manifold reduction to a number of individual-based, discrete stochastic processes that include a number of well-known examples from applications. In \cite{Kurtz70,Kurtz71,Kurtz78,Kurtz81}, Kurtz introduced and studied what he called density dependent population processes.  While his original motivation was chemical reaction networks, the class also includes many examples of interest in biology and epidemiology.

A sequence of Markov processes $\bm{x}_{n}(t)$ is a \textit{density dependent population process} if $\bm{x}_{n}$ takes values in $\frac{1}{n} \mathbb{Z}^{d}$, and, if $q^{(n)}_{\bm{x},\bm{y}}$ is the jump rate between $\bm{x},\bm{y} \in \frac{1}{n} \mathbb{Z}^{d}$, then
\begin{equation*}
q^{(n)}_{\bm{x},\bm{y}} = n \lambda_{n(\bm{y}-\bm{x})}(\bm{x})
\end{equation*}
for some non-negative function $\lambda_{\bm{l}}(\bm{x})$ on $\mathbb{R}^{d}$, where $\bm{l} = n(\bm{y}-\bm{x}) \in \mathbb{Z}^{d}$.  More generally, one can consider the case of functions
 $\lambda^{(n)}_{\bm{l}}(\bm{x})$  that depend on $n$, provided $\lambda^{(n)}_{\bm{l}}(\bm{x})$ converges to a limit $\lambda_{\bm{l}}(\bm{x})$ sufficiently quickly as $n \to \infty$; see \cite{Pollett1990}.

The parameter $n$ corresponds to the ``system size'' in \cite{vanKampen92}, and can be interpreted differently according to the context, as \textit{e.g.} total population size, area, or volume.  For example, consider the stochastic logistic process $X_{n}(t)$ with birth and death rates
\[
Q^{(n)}_{X,X+1} = \beta X\left(1 - \frac{X}{n}\right) \qquad Q^{(n)}_{X,X-1} = \delta X.
\]
Here, $n$ plays the role of the carrying capacity in the deterministic logistic equation, \textit{i.e.} the number of individuals the environment can support: individuals have an intrinsic per-capita birth rate $\beta$, but the offspring will only survive if it arrives in an unoccupied spot in the habitat.  Nondimensionalising, we might consider instead the process $x_{n}(t) = \frac{1}{n} X_{n}(t)$,
with rates
\[
q^{(n)}_{x,x+\frac{1}{n}} = n \beta x(1-x) \qquad q^{(n)}_{x,x-\frac{1}{n}} = n \delta x.
\]
The latter is an example of a density-dependent population process, with
\[
	\lambda_{1}(x) = \beta x (1-x) \qquad \lambda_{-1}(x) = \delta x.
\]

In \cite{Kurtz70}, Kurtz shows that provided
\[
\sum_{\bm{l} \in \mathbb{Z}^{d}} \|\bm{l}\| \sup_{\bm{x} \in \mathcal{K}} \lambda_{\bm{l}}(\bm{x}) < \infty
\]
for all closed and bounded sets $\mathcal{K}$, then if
\[
\bm{f}(\bm{x}) =  \sum_{\bm{l} \in \mathbb{Z}^{d}} \bm{l} \lambda_{\bm{l}}(\bm{x})
\]
is differentiable and $\bm{x}_{n}(t) \to \bm{x}_{0}$, then for any fixed $T > 0$,
\[
\lim_{n \to \infty} \sup_{t \leq T}	|\bm{x}_{n}(t)- \bm{x}(t)| = 0,
\]
where $\bm{x}(t)$ is the solution of $\frac{d\bm{x}}{dt} = \bm{f}(\bm{x})$ with $\bm{x}(0) = \bm{x}_{0}$.

If one assumes that $\lambda_{\bm{l}}(\bm{x})$ is non-zero for only finitely many transitions, say $\bm{l}_{1},\ldots,\bm{l}_{s}$, then, letting $\bm{G}(\bm{x})$ be the matrix with $i$\textsuperscript{th} column
$\bm{l}_{i} \sqrt{\lambda_{\bm{l}_{i}}(\bm{x})}$, $\bm{\eta}(t)$ be an $s$-dimensional It\^o white noise,
and $\bm{z}_{n}(t)$ be the solution of
\[
\frac{d\bm{z}_{n}}{dt} = \bm{f}(\bm{z}_{n}) + \frac{1}{\sqrt{n}} \bm{G}(\bm{z}_{n}) \bm{\eta}(t),
\]
then for any fixed $T > 0$, there exists a constant $C_{T}$ such that
\[
\lim_{n \to \infty} \mathbb{P}\left(\sup_{t \leq T} |\bm{x}_{n}(t) - \bm{z}_{n}(t)|
> \frac{C_{T} \log{n}}{n}\right) = 0.
\]

In our current setting, if $\bm{f}(\bm{x})$ is twice continuously differentiable and once again has a globally attractive $m$-dimensional manifold of equilibria $\Gamma$, then the process $\bm{z}_{n}(t) = \bm{x}_{n}(nt)$ satisfies the conditions of \cite{Katzenberger91}, so that as $n \to \infty$, $\bm{z}_{n}(t)$ converges to a diffusion $\bm{z}(t)$ satisfying equation \eqref{KSDE} for $\bm{f}$ and $\bm{G}$ defined as above.  This result was applied to the study population genetic and epidemiological models in \cite{Parsons2010,Parsons2012}.

\section{Local representations of one-dimensional manifolds}\label{Theta}

In this section, we will discuss how one may obtain a parameterisation $\bm{\gamma}$ of a one-dimensional slow manifold $\Gamma$ and compute the quadratic expansion of the flow field (i.e. the quantities $\bm{v}$ and $\bm{\Theta}$) in the neighbourhood of a point $\bm{x} \in \Gamma$.

%Write $\bm{J}(z)$ for the Jacobian of $\bm{f}$ at the point $\bm{x}=\bm{\gamma}(z)$, that is,
%\begin{equation}
%J_{ij}(z)=\frac{\partial f_i}{\partial x_j}\Bigg|_{\bm{x}=\bm{\gamma}(z)}\,.
%\end{equation}
We start by fixing a basis of generalised eigenvectors of the Jacobian at $\bm{x}_{0}$, $\bm{J}(\bm{x}_{0})$, say $\bm{w}_{1},\ldots,\bm{w}_{n}$, and letting $\bm{W}$ be the corresponding change of basis matrix with the $\bm{w}_{i}$ as columns. Let $\bm{w}_{1}$ to be the eigenvector corresponding to the eigenvalue 0 (which we take to be unique up to scalar multiplication). %This vector is tangent to $\Gamma$ at $\bm{x}=\bm{\gamma}(z)$, hence
Then,
\begin{equation}\label{WJW}
 	\bm{W}^{-1}\bm{J}(\bm{x}_{0})\bm{W} = \begin{bmatrix} 0 & \\ & \bm{J}_{\bm{2}} \end{bmatrix},
\end{equation}
where $\bm{J}_{2}$ is a block-diagonal matrix, with each block acting invariantly on one of the eigenspaces corresponding to the non-zero eigenvalues.

We introduce a new coordinate system
\[
	\bm{z} = \bm{W}^{-1} \left(\bm{x} - \bm{x}_{0}\right).
\]
In this coordinate system, we will construct a parameterisation $\bm{\gamma}(z_{1})$ of
$\Gamma$ such that that $\bm{x}_{0} =\bm{\gamma}(0)$.

In the new coordinate system, the dynamics are then given by $\frac{d\bm{z}}{dt} = \hat{\bm{f}}(\bm{z})$, where
\[
	\hat{\bm{f}}(\bm{z})  = \bm{W}^{-1}\bm{f}\left(\bm{x}_{0} + \bm{W}\bm{z}\right)
\]
(thus, the Jacobian of $\hat{\bm{f}}$ at $\bm{0}$, say $\hat{\bm{J}}$, is $\bm{W}^{-1}\bm{J}(\bm{x}_{0})\bm{W}$).  Setting $\bm{z}_{\bm{2}}=(z_2\,\ldots,z_d)$, we may write this as
\begin{equation}\label{NormalODE}
\begin{aligned}
	\frac{dz_1}{dt} &= \varphi_1(z_1,\bm{z}_{\bm{2}})\\
	\frac{d\bm{z}_{\bm{2}}}{dt} &= \bm{J}_{\bm{2}}\bm{z}_{\bm{2}} + \bm{\varphi}_{\bm{2}}(z,\bm{z}_{\bm{2}})\,.
\end{aligned}
\end{equation}
where $\bm{\varphi}_{\bm{2}}(z_1,\bm{z}_{\bm{2}}) = (\varphi_{2},\ldots,\varphi_{d})$ is quadratic. We may thus Taylor expand $\varphi_{i}(\bm{z})$ about $\bm{0}$ as
\[
	\varphi_{i}(\bm{z}) = \sum_{j, k = 1}^{d} c_{ijk} z_{j} z_{k} + O\left(|\bm{z}|^{3}\right).
\]

Computing $\bm{\gamma}$ or $\bm{\Theta}$, is essentially the task of characterizing the centre and stable manifolds at $\bm{x}_{0}$ respectively. The centre manifold theorem (we follow the treatment in \cite{Glendinning94}) tells us that at $\bm{x}_{0}$ the centre manifold is tangent to $\bm{w}_{1} $, whereas the stable manifold is tangent to the space spanned by $\bm{w}_{2},\ldots,\bm{w}_{d}$.  Moreover, we may locally represent each manifold as the graph of a function over the tangent space.  In particular, in the new coordinate system, there exists a function
\[
	\bm{\gamma}_{\bm{2}}(z_{1}) = (\gamma_{2}(z_{1}),\ldots,\gamma_{d}(z_{1}))
\]
such that $\bm{\gamma}(z) = (z,\bm{\gamma}_{\bm{2}}(z))^T$ is a point on $\Gamma$ for all $z_{1}$ sufficiently close to 0, and a function $\vartheta(\bm{z}_{\bm{2}})$ such that $(\bm{z}_{\bm{2}},\vartheta(\bm{z}_{\bm{2}}))$ is a point in the stable manifold near $\bm{x}_{0}$ for $\bm{z}_{\bm{2}}$ sufficiently close to $\bm{0}$.  We will demonstrate the calculation of $\bm{\gamma}_{\bm{2}}(z_{1})$ below; the calculation of $\vartheta(\bm{z}_{\bm{2}})$ is similar, so we will simply give the result.  Finally, we will show how one obtains $\bm{\Theta}$ from $\vartheta(\bm{z}_{\bm{2}})$.

%The stable (fast) manifold around $\bm{z}$ is (in the local coordinate system)
%\begin{equation}
%\pi^{-1}(\bm{\gamma}(z))=\left\{ \left(\begin{array}{c}z_1\\\bm{z}_{\bm{2}}\end{array}\right)\,:\, z_1=\vartheta(\bm{z}_{\bm{2}})\right\}\,.
%\end{equation}
%So the stable manifold is the graph of the function $\vartheta$.

To begin, we observe that in our new coordinate system $\bm{x}_{0}$ is the origin and $\Gamma$ is tangent to the $z_{1}$ axis (\textit{i.e.}\, the span of $\bm{w}_{1}$), so we must have $\bm{\gamma}_{\bm{2}}'(0) = \frac{d\bm{\gamma}_{\bm{2}}}{dz_{1}} = \bm{0}$.  We thus look for $\bm{\gamma}_{\bm{2}}(z_{1})$ of the form
\[
	\gamma_{i}(z_{1}) = a_{i} z_{1}^{2} + O\left(z_{1}^{3}\right).
\]
(as we shall only be interested in the first and second order derivatives of $\bm{\gamma}$ at $\bm{x}_{0}$ -- \textit{i.e.}\, at $z_{1} = 0$ -- this is adequate for our purposes).

Substituting into \eqref{NormalODE}, for points on $\Gamma$ we have
\begin{equation}
\begin{aligned}
	\frac{dz_{1}}{dt} &= \varphi_{1}(z_{1},\bm{\gamma}_{\bm{2}}(z_{1}))\\
	\frac{d}{dt}\bm{\gamma}_{\bm{2}}(z_{1}) &= \bm{J}_{\bm{2}}\bm{\gamma}_{\bm{2}}(z_{1}) + \bm{\varphi}_{\bm{2}}(z_{1},\bm{\gamma}_{\bm{2}}(z_{1})),
\end{aligned}
\end{equation}
or, expanding the latter using the chain rule,
\[
	\varphi_{1}(z_{1},\bm{\gamma}_{\bm{2}}(z_{1})) \frac{d\bm{\gamma}_{\bm{2}}}{dz_{1}}
	= \bm{J}_{\bm{2}}\bm{\gamma}_{\bm{2}}(z_{1}) + \bm{\varphi}_{\bm{2}}(z_{1},\bm{\gamma}_{\bm{2}}(z_{1})).
\]
Substituting our series expressions for the $\varphi_{i}$ and $h_{i}$, to lowest order this gives us
\[
	2 c_{i11} a_{i} z_{1}^{3} + O\left(z_{1}^{4}\right)
	= \left(\sum_{j=2}^{d} \hat{J}_{ij} a_{j} + c_{i11}\right) z_{1}^{2} + O\left(z_{1}^{3}\right)
\]
\textit{i.e.}\, we may obtain the quantities $a_{i}$, $i=2,\ldots,d$ by solving the system of equations
\[
	\sum_{j=2}^{d} \hat{J}_{ij} a_{j} = - c_{i11},\quad i = 2,\ldots,d.
\]
Noting that $(a_{2},\ldots,a_{d})^{T} = \frac{1}{2} \frac{d^{2}\bm{\gamma}_{\bm{2}}}{dz_{1}^{2}}(0)$ whereas
$(c_{111},\ldots,c_{d11})^{T} = \frac{1}{2} \frac{\partial^{2}\bm{\varphi}_{\bm{2}}}{\partial z_{1}^{2}}(\bm{0})$, we can solve the previous equation as
\[
	\frac{d^{2}\bm{\gamma}_{\bm{2}}}{dz_{1}^{2}}(0)
		= -\hat{\bm{J}}^{-1} \frac{\partial^{2}\bm{\varphi}_{\bm{2}}}{\partial z_{1}^{2}}(\bm{0}).
%		= - \frac{\partial^{2}(\hat{\bm{J}}^{-1}\bm{\varphi}_{\bm{2}})}{\partial z_{1}^{2}}(\bm{0}).
\]
To return to our original functions as expressed in the original coordinate system, we first observe that for $i = 2,\ldots,d$,
\[
	 \frac{d^{2}\varphi_{i}}{dz_{1}^{2}}(\bm{0})
	 =  \frac{\partial^{2}\hat{f}_{i}}{\partial z_{1}^{2}}(\bm{0}),
\]
whereas
\[
	\frac{\partial^{2}\hat{\bm{f}}}{\partial z_{1}^{2}}(\bm{0})
	= \bm{W}^{-1} \sum_{j,k =1}^{d} \frac{\partial^{2} \bm{f}}{\partial x_{j} \partial x_{k}}(\bm{x}_{0}) W_{j1} W_{k1}.
\]
In particular, recalling \eqref{WJW}, we see that for $i = 2,\ldots,d$, $\frac{d^{2}\gamma_{i}}{dz_{1}^{2}}(0)$
agrees with the $i$\textsuperscript{th} entry of
\[
	- \bm{W}^{-1}\bm{J}^{+} \sum_{j,k =1}^{d}
		\frac{\partial^{2} \bm{f}}{\partial x_{j} \partial x_{k}}(\bm{x}_{0}) W_{j1} W_{k1},
\]
\textit{i.e.}
\[
	- \left[\bm{W}^{T}
	\frac{\partial^{2} \left[\bm{W}^{-1}\bm{J}^{+} \bm{f}\right]_{i}}{\partial \bm{x}^{2}}(\bm{x}_{0})
		\bm{W}\right]_{11},
\]
where, as before, $\bm{J}^{+}$ is the pseudo-inverse of $\bm{J}$, which is defined by $\bm{J}^{-1}$ on the image of $\bm{J}$ and is $\bm{0}$ on the kernel of $\bm{J}$.  

Thus,
\begin{gather*}
	\bm{\gamma}(0) = \bm{x}_{0}\\
	\bm{\gamma}'(0) = \bm{w}_{1}\\
	\bm{\gamma}''(0) = \sum_{i=2}^{d} \frac{d^{2}\gamma_{i}}{dz_{1}^{2}}(0) \bm{w}_{i}
	= - \sum_{i=2}^{d}  \left[\bm{W}^{T} \frac{\partial^{2}
	\left[\bm{W}^{-1}\bm{J}^{+} \bm{f}\right]_{i}}{\partial \bm{x}^{2}}(\bm{x}_{0})\bm{W}\right]_{11}\bm{w}_{i} 
	%= - \sum_{i=2}^{d} \frac{\partial^{2}[\hat{\bm{J}}^{-1}\bm{\varphi}_{\bm{2}}]_{i}}{\partial z_{1}^{2}}(\bm{0})\bm{w}_{i} , 
\end{gather*}
and
\[
	\bm{\gamma}(z_1) = \bm{x}_{0} + z_{1}\bm{w}_{1}  
	- \frac{1}{2} z_{1}^{2}  \sum_{i=2}^{d}  \left[\bm{W}^{T} \frac{\partial^{2}
	\left[\bm{W}^{-1}\bm{J}^{+} \bm{f}\right]_{i}}{\partial \bm{x}^{2}}(\bm{x}_{0})\bm{W}\right]_{11}\bm{w}_{i} 
	+ O\left(z_{1}^{3}\right)
\]
is the desired parametrization of $\Gamma$ in the $\bm{z}$ coodinates.

Proceeding similarly, we find that
\[
	\vartheta(\bm{z}_{\bm{2}}) = \bm{z}_{\bm{2}}^{T} (\hat{\bm{J}}^{T})^{-1}
		\frac{\partial^{2} \hat{f}_{1}}{\partial \bm{z}_{\bm{2}}^{2}}(\bm{0}) \bm{z}_{\bm{2}}
\]
and $(\hat{\bm{J}}^{T})^{-1} \frac{\partial^{2} \hat{f}_{1}}{\partial \bm{z}_{\bm{2}}^{2}}(\bm{0})$ and
\[
	\bm{W}^{T}  \left(\bm{J}^{T}\right)^{+}
		\frac{\partial^{2} \left[\bm{W}^{-1}\bm{f}\right]_{1}}{\partial \bm{x}^{2}}(\bm{x}_{0})
		\bm{W}
\]
have equal $j,k$\textsuperscript{th} entry for all $j,k=2,\ldots,d$ (the first row of the latter is zero, but the first column need not be).  Thus, if we set
\[
	\bm{\Theta(\bm{x}_{0})} = \bm{P}(\bm{x}_{0}) \bm{W}^{T}  \left(\bm{J}^{T}\right)^{+}
		\frac{\partial^{2} \left[\bm{W}^{-1}\bm{f}\right]_{1}}{\partial \bm{x}^{2}}(\bm{x}_{0})
		\bm{W},
\]
then the stable manifold at $\bm{x}_{0}$ is thus the set of all points $\bm{z}$ such that
\[
	z_{1} = \bm{z}_{\bm{2}}^{T} \bm{\Theta(\bm{x}_{0})} \bm{z}_{\bm{2}}.
\]
Now, if we choose $\bm{v}(\bm{x}_{0})$ so that
\[
	\bm{v}(\bm{x}_{0})^{T} \bm{w}_{i} = \begin{cases}
		1 & \text{if $i =1$, and}\\
		0 & \text{otherwise,}
	\end{cases}
\]
then for a point $\bm{x} = \bm{x}_{0} + \Delta \bm{x}$, $z_{1} = \bm{v}(\bm{x}_{0})^{T}  \Delta \bm{x}$,
whereas
\[
	\bm{z}_{\bm{2}} = \Delta \bm{x} - \left(\bm{v}(\bm{x}_{0})^{T} \Delta \bm{x}\right) \bm{w}_{1},
\]
so that $\bm{x}$ is in the stable manifold at  $\bm{x}_{0}$ (to lowest order in $\Delta \bm{x}$) provided
\[
	\bm{v}(\bm{x}_{0})^{T}  \Delta \bm{x}
	- \left(\Delta \bm{x} - \left(\bm{v}(\bm{x}_{0})^{T} \Delta \bm{x}\right) \bm{w}_{1}\right)^{T}
	\Theta(\bm{x}_{0})
	\left(\Delta \bm{x} - \left(\bm{v}(\bm{x}_{0})^{T} \Delta \bm{x}\right) \bm{w}_{1}\right) = 0,
\]
or, rearranging,
\[
	\bm{v}(\bm{x}_{0})^{T}  \Delta \bm{x}
	- \Delta \bm{x}^{T}\left(\bm{I} - \bm{w}_{1}\bm{v}(\bm{x}_{0})^{T}\right)^{T}
	\Theta(\bm{x}_{0}) \left(\bm{I} - \bm{w}_{1}\bm{v}(\bm{x}_{0})^{T}\right) \Delta \bm{x}
	= 0.
\]

\section{Derivation of general case}
\label{Q}

First we examine the projection matrix $\bm{P}$. Consider the outer system
\begin{equation}
\frac{d\bm{\xi}}{dt}=\bm{f}(\bm{\xi})\,,\quad \bm{\xi}(0)=\bm{x}\,,
\end{equation}
where $\bm{x}$ lies close to a point $\bm{z}$ on the manifold. Varying the initial conditions yields
\begin{equation}\label{odedxidx}
\frac{d}{dt}\frac{\partial\xi_i}{\partial x_j}=\frac{\partial}{\partial x_j}f_i(\bm{\xi})=\sum_k \frac{\partial\xi_k}{\partial x_j}\frac{\partial }{\partial \xi_k}f_i(\bm{\xi})=\sum_kJ_{ik}(\bm{\xi})\frac{\partial\xi_k}{\partial x_j}
\end{equation}
\textit{i.e.}
\[
	\frac{d}{dt}\frac{\partial \bm{\xi}}{\partial \bm{x}}
	=  \bm{J}(\bm{\xi}(t))\frac{\partial \bm{\xi}}{\partial \bm{x}}
\]
where $\bm{J}$ is the Jacobian matrix of $\bm{f}$.  Now, since $\bm{\xi}(0)=\bm{x}$,
\begin{equation*}
	\frac{\partial \bm{\xi}}{\partial \bm{x}}(0,\bm{x}) = \bm{I}
\end{equation*}
and thus this variational equation has solution
\begin{equation*}
	\frac{\partial\xi_i}{\partial x_j} = \bm{\Pi}(0,t)
\end{equation*}
where $\bm{\Pi}(s,t)$ is the fundamental matrix solving
\[
	\frac{d}{dt} \bm{\Pi}(s,t) = \bm{J}(\bm{\xi}(\bm{x},t))\bm{\Pi}(s,t), \quad \bm{\Pi}(s,s) = I.
\]
When $\bm{x}$ is taken to be $\bm{z} \in \Gamma$, since $\bm{\xi}(\bm{z},t) = \bm{z}$ for all $\bm{z} \in \Gamma$, we have
\[
	\bm{\Pi}(s,t)  = e^{(t-s)\bm{J}(\bm{z})}.
\]
so, in this case, $\frac{\partial \bm{\pi}}{\partial \bm{x}}(\bm{z},t) = e^{t \bm{J}(\bm{z})}$
(\textit{i.e.} informally, $\frac{d}{dt}\frac{\partial \bm{\xi}}{\partial \bm{x}} \approx  \bm{J}(\bm{z})\frac{\partial \bm{\xi}}{\partial \bm{x}}$. Under this approximation the equation is linear and admits the solution $\frac{\partial\bm{\xi}}{\partial \bm{x}}=e^{t\bm{J}(\bm{z})}$).

From the definitions (\ref{defpi}) and (\ref{defPQ}) we recover $\bm{P}$ by taking the limit of large $t$,
\begin{equation}
\bm{P}(\bm{z})=\lim_{t\to\infty}e^{t\bm{J}(\bm{z})}\,.
\end{equation}
To compute the limit we consider the action of $e^{t\bm{J}(\bm{z})}$ on an eigenvector of the Jacobian\footnote{To simplify the discussion we assume that $\bm{J}(\bm{z})$ is diagonalizable and that its kernel contains only the tangent plane to the manifold. Neither assumption is necessary. }. If $\bm{u}_i$ is tangent to the manifold then the corresponding eigenvalue $\lambda_i$ is zero and so $e^{t\lambda_i}=1$ and $\bm{P}(\bm{z})$ leaves $\bm{u}_i$ unchanged. Alternatively, if $\bm{u}_i$ corresponds to a direction of fast collapse then its eigenvalue is negative and $e^{t\lambda_i}\to0$, so $\bm{u}_i$ is annihilated by $\bm{P}(\bm{z})$.

Let $\bm{U}=(\bm{u}_1,\ldots,\bm{u}_m)$ be a basis of the tangent plane to the manifold at $\bm{z}$ (the slow subspace) and let $\bm{V}=(\bm{v}_1,\ldots,\bm{v}_{m})$ a basis of the orthogonal complement of the fast subspace. Then we may write
\begin{equation}
\bm{P}(\bm{z})=\bm{U}(\bm{V}^T\bm{U})^{-1}\bm{V}^T\,.
\label{Puv}
\end{equation}
In the above we assumed that the tangent plane to the manifold was precisely the kernel of the Jacobian, in which case $\bm{U}$ would be the first $m$ columns of the right eigenvector matrix, and $\bm{V}^T$ the bottom $d-m$ rows of the left eigenvector matrix. This may not hold if the manifold is not hyperbolic (for example if $\bm{f}$ has a component like $-x_i^3$, which is stable but not linearly so), however, equation (\ref{Puv}) remains true for all flow fields, provided we somehow have access to bases $\bm{U}$ and $\bm{V}$.

Let us move on to calculate $\bm{Q}$. We start by obtaining some simple identities:  first note that by the definition of $\bm{\pi}$, we have
$f_{i}(\bm{\pi}(\bm{x})) = 0$ for all $\bm{x}$.  Differentiating this, we obtain
\begin{equation}\label{A1}
\sum_{m}\frac{\partial  f_{i}}{\partial x_{m}}(\bm{\pi}(\bm{x})) \frac{\partial \pi_m}{\partial x_j} = 0,
\end{equation}
or, in matrix form, $\bm{J}(\bm{\pi}(\bm{x})) \frac{\partial \bm{\pi}}{\partial \bm{x}} = \bm{0}$.
Replacing $\bm{x}$ by $\bm{z} \in \Gamma$, and recalling that $ \frac{\partial \bm{\pi}}{\partial \bm{x}}(\bm{z}) =\bm{P}(\bm{z})$, we have
\begin{equation*}
\bm{J}(\bm{z})\bm{P}(\bm{z}) = \bm{0},
\end{equation*}
 \textit{i.e.}\, $\bm{J}(\bm{z})$ annihilates all the slow directions, as we have already observed.  Differentiating \eqref{A1}, we obtain
\begin{equation*}
\sum_{m,n}\frac{\partial^{2} f_{i}}{\partial x_{m}\partial x_{n}}(\bm{\pi}(\bm{x})) \frac{\partial \pi_m}{\partial x_j} \frac{\partial \pi_n}{\partial x_k} + \sum_{m}\frac{\partial  f_{i}}{\partial x_{m}}(\bm{\pi}(\bm{x})) \frac{\partial^{2} \pi_m}{\partial x_j \partial x_k} = 0,
\end{equation*}
which we can write in vector form as
\begin{equation}\label{A2}
\bm{\mathcal{H}}_{jk}\left(\frac{\partial \bm{\pi}}{\partial \bm{x}}\right) + \bm{J}(\bm{\pi}(\bm{x})) \frac{\partial^{2} \bm{\pi}}{\partial x_j \partial x_k} = 0,
\end{equation}
where, for any $n \times n$-matrix $A$, $\bm{\mathcal{H}}_{jk}(\bm{A})$ is the vector with $i$\textsuperscript{th} entry
\begin{equation*}
\mathcal{H}_{ijk}(\bm{A}) =  \bm{e}_{j}^{T}\bm{A}^{T} \frac{\partial^{2} f_{i}}{\partial \bm{x}^{2}}\bm{A}\bm{e}_k,
\end{equation*}
where $\bm{e}_j$ is the $j$\textsuperscript{th} standard basis vector, and we have written
$\frac{\partial^{2} f_{i}}{\partial \bm{x}^{2}}$ for the Hessian matrix with $j,k$\textsuperscript{th} entry $\frac{\partial^{2} f_{i}}{\partial x_j \partial x_k}$. \textit{i.e.}, since $\frac{\partial \bm{\pi}}{\partial x_j} = \frac{\partial \bm{\pi}}{\partial \bm{x}} \bm{e}_j$,
\begin{equation*}
\mathcal{H}_{ijk}\left(\frac{\partial \bm{\pi}}{\partial \bm{x}}\right)
= \left(\frac{\partial \bm{\pi}}{\partial x_j}\right)^T \frac{\partial^{2} f_{i}}{\partial \bm{x}^{2}}\frac{\partial \bm{\pi}}{\partial x_k}
= \sum_{m,n}\frac{\partial^{2} f_{i}}{\partial x_{m}\partial x_{n}}(\bm{\pi}(\bm{x})) \frac{\partial \pi_m}{\partial x_j} \frac{\partial \pi_m}{\partial x_k}.
\end{equation*}

Now, recalling that at $\bm{z} \in \Gamma$, $\bm{\pi}(\bm{z}) = \bm{z}$, $\frac{\partial \bm{\pi}}{\partial \bm{x}} = \bm{P}(\bm{z})$, and $\frac{\partial^{2} \pi_i}{\partial x_j \partial x_k}(\bm{z}) = Q_{ijk}(\bm{z})$, we can write \eqref{A2}  as
\begin{equation}\label{A3}
\bm{J}(\bm{z}) \bm{Q}_{jk}(\bm{z}) = - \bm{\mathcal{H}}_{jk}(\bm{P}(\bm{z})),
\end{equation}
where we continue with the convention that $\bm{Q}_{jk}(\bm{z})$ is the vector with $i$\textsuperscript{th} entry $Q_{ijk}(\bm{z})$.  Applying $\bm{P}(\bm{z})$ to both sides of \eqref{A3} gives
\begin{equation*}
\bm{P}(\bm{z})\bm{\mathcal{H}}_{jk}(\bm{P}(\bm{z})) = \bm{0}
\end{equation*}
so we see $\bm{\mathcal{H}}_{jk}(\bm{P}(\bm{z}))$ is entirely contained in the eigenspace of fast directions.  Notice that restricted to the fast subspace, $\bm{J}(\bm{z})$ is a full-rank operator, so that, \emph{regarded as an operator on the fast subspace}, \eqref{A3} has a unique solution, which we will write as
\begin{equation*}
- \bm{J}(\bm{z})^{+} \bm{\mathcal{H}}_{jk}(\bm{P}(\bm{z})).
\end{equation*}
where we recall that $\bm{J}(\bm{z})^{+}$ is the pseudo-inverse of $\bm{J}(\bm{z})$,  which acts as the inverse of $\bm{J}(\bm{z})$ when restricted to the fast directions and which annihilates all vectors in the slow directions.
%($\bm{J}(\bm{z})^{+}$ is an example of a pseudo-inverse: $\bm{J}(\bm{z})$ is not invertible, but we can compute the matrix $\bm{J}(\bm{z})^{+}$ which acts as the inverse of $\bm{J}(\bm{z})$ when restricted to the fast directions and which annihilates all vectors in the slow directions.)  

However, regarded as an equation on all of $\mathbb{R}^{d}$, the solution to \eqref{A3} is not unique, but rather takes the form
\begin{equation*}
\bm{Q}_{jk}(\bm{z}) = - \bm{J}(\bm{z})^{+} \bm{\mathcal{H}}_{jk}(\bm{P}(\bm{z})) + \bm{S}_{jk}(\bm{z})
\end{equation*}
for some vector $\bm{S}_{jk}(\bm{z})$ in the slow directions.

To obtain $\bm{S}_{jk}(\bm{z})$, we proceed as we did to obtain $\bm{P}(\bm{z})$, differentiating \eqref{odedxidx} to obtain
\begin{equation*}
\frac{d}{dt}\frac{\partial^{2}\xi_i}{\partial x_j \partial x_k}
= \sum_l \bm{J}_{il}(\bm{\xi})\frac{\partial^{2}\xi_l}{\partial x_j \partial x_k}
+\sum_{m,n}\frac{\partial^{2} f_{i}}{\partial x_{m}\partial x_{n}}(\bm{\xi})\frac{\partial\xi_m}{\partial x_j}\frac{\partial\xi_n}{\partial x_k}
\end{equation*}
which again write in vector form as
\begin{equation}\label{QPDE}
\frac{d}{dt}\frac{\partial^{2}\bm{\xi}}{\partial x_j \partial x_k}
 = \bm{J}(\bm{\xi}) \frac{\partial^{2}\bm{\xi}}{\partial x_j \partial x_k} + \bm{\mathcal{H}}_{jk}\left(\frac{\partial \bm{\xi}}{\partial \bm{x}}\right).
\end{equation}
This may be formally solved by Duhamel's principle to give
\begin{equation*}
\frac{\partial^{2}\bm{\xi}}{\partial x_j \partial x_k} =
\int_{0}^{t} \bm{\Pi}(s,t)%e^{\int_{s}^{t} \bm{J}(\bm{\xi}(\bm{x},u))\, du}
\bm{\mathcal{H}}_{jk}\left(\frac{\partial \bm{\xi}}{\partial \bm{x}}(x,s)\right)\, ds
\end{equation*}
where $\bm{\Pi}(s,t)$ is the fundamental matrix from above.

As before, when $\bm{x}$ is taken to be a point $\bm{z} \in \Gamma$, since $\bm{\xi}(\bm{z},t) = \bm{z}$ for all $\bm{z} \in \Gamma$, we have $\bm{\Pi}(s,t)  = e^{(t-s)\bm{J}(\bm{z})}$ and $\frac{\partial \bm{\pi}}{\partial \bm{x}}(\bm{z},t) = e^{t \bm{J}(\bm{z})}$, and the solution to \eqref{QPDE} simplifies to
\begin{equation*}
\int_{0}^{t} e^{(t-s) \bm{J}(\bm{z})} \bm{\mathcal{H}}_{jk}\left(e^{s\bm{J}(\bm{z})}\right)\, ds.
\end{equation*}
Thus,
\begin{equation*}
\bm{Q}_{jk}(\bm{z})
= \lim_{t \to \infty} \frac{\partial^{2}\bm{\xi}}{\partial x_j \partial x_k}
= \lim_{t \to \infty} \int_{0}^{t} e^{(t-s) \bm{J}(\bm{z})} \bm{\mathcal{H}}_{jk}\left(e^{s\bm{J}(\bm{z})}\right)\, ds.
\end{equation*}

Now, $e^{t\bm{J}(\bm{z})} \to\bm{P}(\bm{z})$ as $t \to \infty$, and
\begin{equation*}
\lim_{t \to \infty} \bm{\mathcal{H}}_{jk}\left(e^{t\bm{J}(\bm{z})}\right)
= \bm{\mathcal{H}}_{jk}(\bm{P}(\bm{z})),
\end{equation*}
both of which are non-zero, so it is not immediately obvious that the integral above converges.  However, the information obtained above allows us to resolve these issues.  We start by observing that
\begin{equation*}
\bm{S}_{jk}(\bm{z}) =\bm{P}(\bm{z})\bm{Q}_{jk}(\bm{z}) =  \lim_{t \to \infty} \int_{0}^{t}\bm{P}(\bm{z}) e^{(t-s) \bm{J}(\bm{z})}  \bm{\mathcal{H}}_{jk}\left(e^{s\bm{J}(\bm{z})}\right)\, ds,
\end{equation*}
and, since $e^{(t-s) \bm{J}(\bm{z})}$ acts like the identity matrix on the slow directions, $\bm{P}(\bm{z}) e^{(t-s) \bm{J}(\bm{z})} =\bm{P}(\bm{z})$, so that
\begin{equation*}
\bm{S}_{jk}(\bm{z}) = \lim_{t \to \infty} \int_{0}^{t} \bm{P}(\bm{z}) \bm{\mathcal{H}}_{jk}\left(e^{s\bm{J}(\bm{z})}\right)\, ds =\bm{P}(\bm{z}) \int_{0}^{\infty} \bm{\mathcal{H}}_{jk}\left(e^{s\bm{J}(\bm{z})}\right)\, ds.
\end{equation*}
Moreover, we've already observed that $\bm{P}(\bm{z})\bm{\mathcal{H}}_{jk}(\bm{P}(\bm{z})) = \bm{0}$, so
\begin{equation*}
\bm{S}_{jk}(\bm{z}) =\bm{P}(\bm{z}) \int_{0}^{\infty} \bm{\mathcal{H}}_{jk}\left(e^{s\bm{J}(\bm{z})}\right) - \bm{\mathcal{H}}_{jk}(\bm{P}(\bm{z}))\, ds,
\end{equation*}
and we are left with evaluating the integral
\begin{equation*}
\begin{multlined}
\int_{0}^{\infty} \mathcal{H}_{ijk}\left(e^{s\bm{J}(\bm{z})}\right) - \mathcal{H}_{ijk}(\bm{P}(\bm{z}))\, ds
= \int_{0}^{\infty}  \bm{e}_j^T e^{s\bm{J}(\bm{z})^{T}} \frac{\partial^{2} f_{i}}{\partial \bm{x}^{2}}(\bm{z}) e^{s\bm{J}(\bm{z})} \bm{e}_k
-  \bm{e}_j^T\bm{P}(\bm{z})^T \frac{\partial^{2} f_{i}}{\partial \bm{x}^{2}}(\bm{z})\bm{P}(\bm{z}) \bm{e}_k\\
=  \bm{e}_j^T \left(\int_{0}^{\infty} e^{s\bm{J}(\bm{z})^{T}} \frac{\partial^{2} f_{i}}{\partial \bm{x}^{2}}(\bm{z}) e^{s\bm{J}(\bm{z})} -\bm{P}(\bm{z})^T \frac{\partial^{2} f_{i}}{\partial \bm{x}^{2}}(\bm{z})\bm{P}(\bm{z})\,ds\right) \bm{e}_k.
\end{multlined}
\end{equation*}
Now,
\begin{equation*}
\begin{multlined}
\int_{0}^{\infty} e^{s\bm{J}(\bm{z})^{T}} \frac{\partial^{2} f_{i}}{\partial \bm{x}^{2}}(\bm{z}) e^{s\bm{J}(\bm{z})} -\bm{P}(\bm{z})^T \frac{\partial^{2} f_{i}}{\partial \bm{x}^{2}}(\bm{z})\bm{P}(\bm{z})\,ds\\
= \int_{0}^{\infty} (e^{s\bm{J}(\bm{z})}-\bm{P}(\bm{z}))^{T} \frac{\partial^{2} f_{i}}{\partial \bm{x}^{2}}(\bm{z}) (e^{s\bm{J}(\bm{z})}-\bm{P}(\bm{z}))\,ds\\
+ \int_{0}^{\infty} (e^{s\bm{J}(\bm{z})}-\bm{P}(\bm{z}))^{T} \frac{\partial^{2} f_{i}}{\partial \bm{x}^{2}}(\bm{z})\bm{P}(\bm{z})\,ds
+ \int_{0}^{\infty}\bm{P}(\bm{z})^{T} \frac{\partial^{2} f_{i}}{\partial \bm{x}^{2}}(\bm{z}) (e^{s\bm{J}(\bm{z})}-\bm{P}(\bm{z}))\,ds,
\end{multlined}
\end{equation*}
and, since $e^{t\bm{J}(\bm{z})}-\bm{P}(\bm{z})$ vanishes on the slow directions, and acts as $e^{t\bm{J}(\bm{z})}$ restricted to the fast directions,
\begin{equation*}
\int_{0}^{\infty} e^{s\bm{J}(\bm{z})}-\bm{P}(\bm{z})\, ds = -\bm{J}(\bm{z})^{+}
\end{equation*}
whereas
\begin{equation}\label{lyapunovsol}
\bm{X}_i(\bm{z}) = \int_{0}^{\infty} (e^{s\bm{J}(\bm{z})}-\bm{P}(\bm{z}))^{T} \frac{\partial^{2} f_{i}}{\partial \bm{x}^{2}}(\bm{z}) (e^{s\bm{J}(\bm{z})}-\bm{P}(\bm{z}))\,ds
\end{equation}
is the unique solution to the Lyapunov equation
\begin{equation*}
\bm{J}(\bm{z})^{T} \bm{X}_i(\bm{z}) + \bm{X}_i(\bm{z}) \bm{J}(\bm{z})
= - \frac{\partial^{2} f_{i}}{\partial \bm{x}^{2}}(\bm{z})
\end{equation*}
in the fast subspace \cite{Bellman60}.  Thus,
\begin{equation*}
\bm{S}_{jk}(\bm{z}) =\bm{P}(\bm{z}) \tilde{\bm{S}}_{jk}(\bm{z}),
\end{equation*}
where
\begin{equation*}
\tilde{S}_{ijk}(\bm{z}) = \bm{e}_j^T \left(\bm{X}_i(\bm{z})
- (\bm{J}(\bm{z})^{+})^{T}\frac{\partial^{2} f_{i}}{\partial \bm{x}^{2}}\bm{P}(\bm{z})
- \bm{P}(\bm{z})^{T} \frac{\partial^{2} f_{i}}{\partial \bm{x}^{2}}\bm{J}(\bm{z})^{+}\right) \bm{e}_k
\end{equation*}
and, finally,
\begin{equation*}
\bm{Q}_{jk}(\bm{z})
= - \bm{J}(\bm{z})^{+} \bm{\mathcal{H}}_{jk}(\bm{P}(\bm{z})) + \bm{P}(\bm{z}) \tilde{\bm{S}}_{jk}(\bm{z}).
\end{equation*}

%\section{Systems with invariants of motion}

%\section{General Riemann manifolds}

%\section{Notation table}
%\begin{center}
%\begin{tabular}{p{2cm}|p{10cm}|p{3cm}}
%Symbol & Meaning & Place of defintion\\\hline
%$d$ & Dimension of space & \\
%$\bm{x}$ & Co-ordinate in $\mathbb{R}^n$, and stochastic variable evolving according to equation (1) & \\
%\end{tabular}
%\end{center}

\bibliography{DRvTS}
\bibliographystyle{apalike}

\end{document}